\documentclass[12pt]{amsart}
\usepackage{amsmath}
\usepackage{amssymb}
\usepackage{amsfonts}
\usepackage{amsthm}
\usepackage{mathrsfs}
\usepackage{color}
\usepackage{graphicx}
\usepackage{caption}
\usepackage[all]{xy}
\theoremstyle{remark}
\newtheorem*{theorem*}{Theorem}
\newtheorem{theorem}{Theorem}[section]
\newtheorem{lemma}[theorem]{Lemma}

\newtheorem{proposition}[theorem]{Proposition}

\usepackage{graphicx}
\graphicspath{{images/}}

\theoremstyle{remark}
\newtheorem{definition}[theorem]{Definition}
\newtheorem{remark}[theorem]{Remark}

\newtheorem{example}[theorem]{Example}

\begin{document}
\title{Modular representations in type $A$ with a two-row nilpotent central character}

\author[G.~Dobrovolska]{Galyna Dobrovolska}
\address{Department of Mathematics, Ariel University, Ariel 40700 Israel}
\email{galdobr@gmail.com}

\author[V.~Nandakumar]{Vinoth Nandakumar}
\address{School of Mathematics and Statistics, University of Sydney, NSW 2006 Australia}
\email{vinoth.nandakumar@sydney.edu.au}

\author[D.~Yang]{David Yang}
\address{Department of Mathematics, Harvard University, One Oxford Street, Cambridge MA 02138}
\email{dyang@math.harvard.edu}
\dedicatory{Dedicated to our friend, Yi Sun}
\maketitle

\begin{abstract} We study the category of representations of $\mathfrak{sl}_{m+2n}$ over a field of characteristic $p$ with $p \gg 0$, whose $p$-character is a nilpotent whose Jordan type is the two-row partition $(m+n,n)$. In a previous paper with Anno, we used Bezrukavnikov-Mirkovic-Rumynin's theory of positive characteristic localization and exotic t-structures to give a geometric parametrization of the simples using annular crossingless matchings. Building on this, here we give combinatorial dimension formulae for the simple objects, and compute the Jordan-H\"{o}lder multiplicities of the simples inside the baby Vermas. We use Cautis-Kamnitzer's geometric categorification of the tangle calculus to study the images of the simple objects under the BMR equivalence. Our results generalize Jantzen's formulae in the subregular nilpotent case (i.e. when $n = 1$), and may be viewed as a positive characteristic analogue of the combinatorial description for Kazhdan-Lusztig polynomials of Grassmannian permutations due to Lascoux and Schutzenberger. \end{abstract}

\section{Introduction}


Representations of algebraic groups and Lie algebras over fields of positive characteristic have been the subject of much study in representation theory. When the characteristic is sufficiently large, dimension formulae for the simple objects can be deduced from Lusztig's conjectures. These are analogous to the Kazhdan-Lusztig conjecture for characters of simple objects in category $\mathcal{O}$, and their proofs also rely on geometric localization techniques. Lusztig's conjecture for algebraic groups gives dimension formulae for the irreducible representations via affine Kazhdan-Lusztig polynomials, and have been established by Andersen-Jantzen-Soergel, \cite{ajs}, building on earlier work of Kashiwara-Tanisaki and Kazhdan-Lusztig. Lusztig's conjectures for Lie algebras study the category of representations with a fixed nilpotent $p$-character, and state the classes of the simple objects match up with canonical bases in Grothendieck groups of Springer fibers. These have been established by Bezrukavnikov and Mirkovic in \cite{bm}, building on their earlier work with Rumynin in \cite{bmr}. 

While the dimensions are in principle known via Lusztig's conjectures, in examples these are computationally difficult. One may ask whether there are any cases where tractable combinatorial formulae exist. In this paper, we show that such a description exists for the irreducible representations of $\mathfrak{sl}_{m+2n}$ whose $p$-character is a nilpotent with Jordan type $(m+n, n)$ (i.e. a two-row partition). We also give combinatorial formulae for the multiplicities of the simples in the baby Vermas in this case.

This is loosely analogous to the principal block of parabolic category $\mathcal{O}^{\mathfrak{p}}$, for the parabolic $\mathfrak{p}$ with Levi $\mathfrak{gl}_m \oplus \mathfrak{gl}_n$ inside $\mathfrak{gl}_{m+n}$. While Kazhdan-Lusztig polynomials are in general difficult to compute, for this block there are explicit combinatorial formulae for the characters of the simple modules (see Lascoux-Schutzenberger, \cite{ls}, and for an exposition see Section 2.4 of Bernstein-Frenkel-Khovanov, \cite{bfk}). Multiplicities of the simples in the Verma objects follow from these character formulae, and are known to be either zero or one. In \cite{bs}, Brundan and Stroppel show that this block of parabolic category $\mathcal{O}$ is equivalent to modules over a variant of Khovanov's arc algebra (defined in their previous paper \cite{bs1}), and use this to reprove these character and multiplicity formulae.

The proof of Lusztig's conjectures in \cite{bm} builds upon Bezrukavnikov-Mirkovic-Rumynin's positive characteristic localization theory from \cite{bmr}, which gives a derived equivalence between categories of modular representations, and categories of coherent sheaves on Springer theoretic varieties. In a previous paper joint with Rina Anno, \cite{an}, we use [BMR] localization theory, combined with Cautis and Kamnitzer's tangle categorification results from \cite{ck}, to study the case where $\mathfrak{g}=\mathfrak{sl}_{m+2n}$ and the $p$-character is a nilpotent with Jordan type $(m+n, n)$. Under the equivalences from \cite{bmr}, the irreducible representations will correspond to complexes of coherent sheaves, which are known as ``simple exotic sheaves''. We show that they are indexed by combinatorial objects that we call ``annular crossingless matchings'', and give a construction of them using the functorial tangle invariants in \cite{ck}.

We use our results from \cite{an} to calculate the dimensions of the irreducible representations. [BMR]-localization theory allows to express the dimensions in terms of the Euler characteristics of the simple exotic sheaves; we compute the latter quantity on the Grothendieck group using Cautis-Kamnitzer's tangle categorification from \cite{ck}. While functorial tangle invariants also appear in Bernstein-Frenkel-Khovanov's construction \cite{bfk} in the category $\mathcal{O}$ set-up, the key difference here is the action of the affine tangle calculus constructed in \cite{an}. 


\subsection{Overview of results} Let $\mathfrak{g} = \mathfrak{sl}_{m+2n}$ be defined over an algebraically closed field $\textbf{k}$ of characteristic $p \gg 0$.
Let $\text{Mod}^{\text{fg}}_{e, \lambda}(U \mathfrak{g})$ be the full category of all finitely generated modules, where the Frobenius center $Z_{\text{Fr}}$ acts with a generalized nilpotent character $e \in \mathfrak{g}^*$, and the Harish-Chandra center $Z_{\text{HC}}$ acts via a fixed central character $\lambda \in \mathfrak{h}^* /\!/ W$
(see Section \ref{modularrep} for more details). Consider the case where the Jordan type of $e$ is a partition with two rows; we refer to these as two-row nilpotents. The main result of this paper is a computation of the dimensions of the simple objects. In order to state our results, we briefly recall the combinatorial set-up from \cite{an}, and the parametrization of the irreducibles there.  

\begin{figure} 
\centering
\captionsetup{justification=centering}
\includegraphics[scale=0.75]{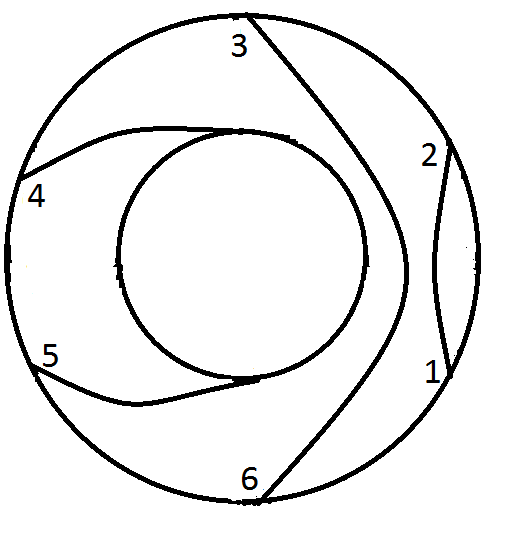}
\caption{}
\label{example}
\end{figure}

\begin{definition} \label{matching} Let $\text{Cross}(m,n)$ be the set of all matchings on an annulus with $m$ unlabelled points on the inner circle, $m+2n$ labelled points on the outer circle (with labels $1, \cdots, m+2n$), such that there are no crossings or loops, and each point on the inner circle is connected to a point on the outer circle. \end{definition}

\begin{example} \label{cross22} See Figure \ref{example} for an example of an element $\alpha \in \text{Cross}(2,2)$. \end{example}

In \cite{an}, we show that the simple objects in $\text{Mod}^{\text{fg}}_{e, \lambda}(U \mathfrak{g})$ can be parametrized by $\text{Cross}(m,n)$. Let us denote the simple objects in $\text{Mod}^{\text{fg}}_{e, 0}(U \mathfrak{g})$ by $M_{\alpha}$ (see Section \ref{[AN]} for more details). Note that this is different from the more general parametrization of the simple objects described in Section \ref{modularrep} below; however the conjecture at the end of Section \ref{conj} below describes how the two are related. Before stating our formula for $\text{dim}(M_{\alpha})$, we introduce some more notation.  

\begin{definition} \label{Calpha} Given $\alpha \in \text{Cross}(m,n)$, define the set $C(\alpha)$ to consist of the following pairs: $$ C(\alpha) = \{ (i, j) \; | \; 1 \leq i < j \leq m+2n; i \text{ and } j \text{ are connected in the outer circle on the diagram } \alpha \} $$ \end{definition}
\begin{definition} Fix the ``dotted line'' to be any line connecting a point lying between $m+2n$ and $1$ in the outer circle, to a point in the inner circle.  Define the disjoint union decomposition $C(\alpha) = C_0(\alpha) \sqcup C_1(\alpha)$ as follows: declare that $(i,j)$ lies in $C_1(\alpha)$ if the dotted line crosses the arc connecting $i$ and $j$, and that $(i,j)$ lies in $C_0(\alpha)$ otherwise. \end{definition}

Let $\mathfrak{h}_{\textbf{k}}$ be the Cartan subalgebra of $\mathfrak{g}$ over $\textbf{k}$. Let $\mathfrak{h}_{\mathbb{Z}}$ be the $\mathbb{Z}$-span of the fundamental weights. Let $\tilde{\mu} = (\tilde{\mu}_1, \cdots, \tilde{\mu}_{m+2n}) \in \mathfrak{h}^*_{\mathbb{Z}}$ lie in the fundamental alcove (see Proposition \ref{dimension} for an explanation of why this restriction is necessary). Let $\mu = (\mu_1, \cdots, \mu_{m+2n}) \in \mathfrak{h}^*_{\textbf{k}}$ be its image under the natural map $\mathfrak{h}^*_{\mathbb{Z}} \rightarrow \mathfrak{h}^*_{\textbf{k}}$. Denote by $T_{0 \rightarrow \mu}: \text{Mod}^{fg}_{e, 0}(\text{U}\mathfrak{g}) \rightarrow \text{Mod}^{fg}_{e, \mu}(\text{U}\mathfrak{g})$ the translation functor. Then $\{ T_{0 \rightarrow \mu} M_{\alpha} \; | \; \alpha \in \text{Cross}(m,n) \} $ constitute the simple objects in $\text{Mod}^{fg}_{e, \mu}(\text{U}\mathfrak{g})$. Our main result is that: 

\begin{theorem}[Main Theorem] \label{main}
\begin{align*} \text{dim}(T_{0 \rightarrow \mu} M_{\alpha}) &= p^{\frac{1}{2}(m+2n)(m+2n-1) - n} \prod_{(i, j) \in C_0(\alpha)} (\tilde{\mu}_i - \tilde{\mu}_j - i +j) \prod_{(i, j) \in C_1(\alpha)} (p-\tilde{\mu}_i + \tilde{\mu}_j + i - j) \end{align*} \end{theorem}

Taking $\mu = 0$, this implies that $$ \text{dim}(M_{\alpha}) = p^{\frac{1}{2}(m+2n)(m+2n-1) - n} \prod_{(i, j) \in C_0(\alpha)} (j - i) \prod_{(i, j) \in C_1(\alpha)} (p + i - j) $$

\begin{example} \label{examp} Let $\alpha \in \text{Cross}(2,2)$ be as above. Then $C_0(\alpha) = \{ (1,2) \}$, $C_1(\alpha) = \{ (3,6) \}$, and: $$ \text{dim}(M_{\alpha}) = p^{13}(p-3), \qquad \text{dim}(T_{0 \rightarrow \mu}M_{\alpha}) = p^{13}(\tilde{\mu}_1 - \tilde{\mu}_2 +1)(p - \tilde{\mu}_3 + \tilde{\mu}_6 - 3 ) $$ \end{example}

We also determine the Jordan-H\"{o}lder multiplicities of the simple objects in the baby Vermas.

\subsection{Contents.} Now we summarize the contents of the paper. 

In Section \ref{prelim}, we review some background material on modular representation theory of Lie algebras, BMR localization theory (cf. \cite{bmr}), and the results from the previous paper \cite{an} constructing the exotic sheaves corresponding to the simple objects. 

In Section \ref{cohsheaves}, we collect some results about the category of coherent sheaves on the Cautis-Kamnitzer varieties $Y_n$ that will be needed in subsequent sections. Most of the results in this section are known to experts, but we state them for the reader's convenience.

In Section \ref{dimensions}, we prove the Main Theorem using the following key fact from \cite{bmr}: if $\Psi_{\alpha}$ is the exotic sheaf corresponding to $M_{\alpha}$, then $\text{dim}(T_{0 \rightarrow \mu} M_{\alpha}) = \chi(\text{Fr}^* \Psi_{\alpha} \otimes \mathcal{O}(\mu + p \rho))$. The latter quantity can be computed once we know the following: the Euler characteristic as a function from the Grothendieck group to $\mathbb{Z}$ (see Lemma \ref{Eulerchar}), the classes of $\Psi_{\alpha}$ in the Grothendieck group (see Lemma \ref{class_easy}), and the effects of Frobenius pullback and tensoring with line bundles as operators on the Grothendieck group (see Lemma \ref{6.2.6} and Lemma \ref{computation}).


In Section \ref{multip}, we give combinatorial formulae for the Jordan-H\"{o}lder multiplicities of simple objects inside the baby Vermas (see Theorem \ref{mult}). This is done by looking at their classes within the Grothendieck group, and computing the transition matrix.


\subsection{Acknowledgements} We would like to thank Roman Bezrukavnikov, Joel Kamnitzer and Sabin Cautis for many helpful conversations. We are also grateful to Rina Anno, Jonathan Brundan, Michael Ehrig, Simon Goodwin, Lewis Topley and Ben Webster for useful comments. The first author would also like to thank the University of Utah (in particular, Peter Trapa) for supporting this research. 

\section{Preliminaries} \label{prelim}

\subsection{Modular representation theory} \label{modularrep}

In this subsection we collect some facts about the representation theory of semisimple Lie algebras over fields of positive characteristic, and refer the reader to Jantzen's expository article \cite{jantzen2} for a detailed treatment. 

Let $G$ be a semisimple, simply connected, algebraic group, with Lie algebra $\mathfrak{g}$, defined over a field $\textbf{k}$ of characteristic $p \gg 0$ (but many of the below results are true under weaker assumptions on $p$; see conditions (H1)-(H3) in B.6 of \cite{jantzen2}). Let $\mathfrak{g} = \mathfrak{n}^- \oplus \mathfrak{h} \oplus \mathfrak{n}^+$ be a triangular decomposition, $\mathfrak{b} = \mathfrak{h} \oplus \mathfrak{n}^+$ and $B \subset G$ be the corresponding Borel. Let $W$ be the associated Weyl group, and $\rho$ the half-sum of all positive roots. Recall that we have the twisted action of $W$ on $\mathfrak{h}^*$: $$w \cdot \lambda = w(\lambda + \rho) - \rho$$

See \cite[\S~B,C]{jantzen2} for a more detailed exposition of the below statements.

\begin{definition} Define the Harish-Chandra center $Z_{\text{HC}}$ to be $Z_{\text{HC}} = (\text{U}\mathfrak{g})^G$. Given an element $x \in \mathfrak{g}$, it is known that there exists a unique $x^{[p]} \in \mathfrak{g}$ such that $x^p - x^{[p]} \in Z(U \mathfrak{g})$. Define the Frobenius center $Z_{\text{Fr}}$ to be the subalgebra generated by $\{ x^p - x^{[p]} \; | \; x \in \mathfrak{g} \}$. When $p \gg 0$, the center of the universal enveloping algebra, $Z(U \mathfrak{g})$, is generated by $Z_{\text{Fr}}$ and $Z_{\text{HC}}$. \end{definition} 

Now we can define the representation categories that we are interested in:

\begin{definition} A character $\lambda \in \mathfrak{h}^*//W$ is integral if it can be expressed as an integral combination of fundamental weights. An integral central character $\lambda \in \mathfrak{h}^*//W$ is said to be regular if its (twisted) $W$-orbit contains $|W|$ elements, and singular otherwise. \end{definition} 

\begin{definition} Suppose $e \in \mathfrak{g}^*$ is a nilpotent character that satisfies $e|_{\mathfrak{b}} = 0$. Let $\text{Mod}^{\text{fg}}_{e}(U \mathfrak{g})$ be the full category of all finitely generated modules, where the Frobenius center $Z_{\text{Fr}}$ acts with a generalized nilpotent character $e \in \mathfrak{g}^*$. Let $\text{Mod}^{\text{fg}}_{e, \lambda}(U \mathfrak{g})$ be the full category of all finitely generated modules, where the Frobenius center $Z_{\text{Fr}}$ acts with a generalized nilpotent character $e \in \mathfrak{g}^*$, and the Harish-Chandra center $Z_{\text{HC}}$ acts via a fixed central character $\lambda \in \mathfrak{h}^* /\!/ W$. We will refer to $\text{Mod}^{\text{fg}}_{e, \lambda}(U \mathfrak{g})$ as a block, and say that the block is singular (or regular) depending on whether $\lambda$ is singular (or regular). \end{definition}

The following decomposition (see Section C of \cite{jantzen2}) justifies the terminology ``block''. $$\text{Mod}^{\text{fg}}_e(U\mathfrak{g}) \simeq \bigoplus_{\lambda \in \mathfrak{h}^*//W} \text{Mod}^{\text{fg}}_{e, \lambda}(U \mathfrak{g})$$

The translation functors $T_{0 \rightarrow \lambda}: \text{Mod}_{e, 0}(U \mathfrak{g}) \rightarrow \text{Mod}_{e, \lambda}(U \mathfrak{g})$ are defined similarly to the category $\mathcal{O}$ case; if $\lambda$ is regular, $T_{0 \rightarrow \lambda}$ is an equivalence (see Jantzen, \cite{jantzenT}, for a detailed discussion).

\begin{definition} Define $U_e(\mathfrak{g})$ to be the quotient of $U\mathfrak{g}$ by the ideal $\langle x^p - x^{[p]} - e(x)^p \; | \; x \in \mathfrak{g} \rangle$. Let $\mu \in \mathfrak{h}^*$ be an integral weight. Define $U_e(\mathfrak{b})$ to be the quotient of $U\mathfrak{b}$ by the ideal $\langle x^p - x^{[p]} - e(x)^p \; | \; x \in \mathfrak{b} \rangle$, and $\textbf{k}_{\mu}$ the one-dimensional $U_e(\mathfrak{b})$-module with highest weight $\mu$ (it is well-defined since $e|_{\mathfrak{b}} = 0$). Then the baby Verma module $Z_e(\mu)$ is defined as: $$Z_e(\mu) = U_e(\mathfrak{g}) \otimes_{U_e(\mathfrak{b})} \textbf{k}_{\mu}$$ \end{definition}

In order to classify the simple modules, we will need to make the assumption that $e$ is in standard Levi form (see Section D1 of \cite{jantzen2} for a definition; this is always the case in type A). The module $Z_e(\mu)$ has a unique maximal submodule, and we denote the head by $L_e(\mu)$.  Unlike the category $\mathcal{O}$ case, it is possible for different baby Verma objects to have the same head. Baby Verma modules share many properties with Verma modules in category $\mathcal{O}$. However, the following proposition illustrates one important difference (see Humphreys, \cite{hum2}, for a proof when $e=0$; and Section C.2 of \cite{jantzen2} for the general case).   

\begin{proposition} If $\mu, \mu' \in W \cdot \lambda$, the Jordan-H\"{o}lder composition multiplicity $[Z_e(\mu): L_e(\mu')]$ does not depend on the choice of $\mu$, i.e. it only depends on $\mu'$. \end{proposition}

\begin{definition} Let $I$ be the corresponding subset of the root system, and let $W(e) := W_I$ be the parabolic Weyl group. Let $W^L(e)$ be a set of maximal-length left coset representatives of $W(e)$ in $W$. \end{definition}

\begin{proposition} Every simple module in $\text{Mod}^{\text{fg}}_{e, \lambda}(U \mathfrak{g})$ has the form $L_e(\mu)$ for some $\mu \in W \cdot \lambda$. Then: 
$$ L_e(\mu) \simeq L_e(w \cdot \mu) \Leftrightarrow w \in W(e) $$
In other words, the irreducible objects in $\text{Mod}^{\text{fg}}_{e, \lambda}(U \mathfrak{g})$ can be parametrized as follows: $\{ L_e(w \cdot \lambda) \; | \; w \in W^L(e) \}$. \end{proposition}

\begin{remark} When we wish to emphasize that the simple module is defined over the field $\textbf{k}$, we will denote it using the notation $L_e(\mu)_{\textbf{k}}$ . \end{remark} 

See Section 3.1 of \cite{nzh} for an exposition of the $\mathfrak{g}=\mathfrak{sl}_2$ example.

\subsection{[BMR] localization theory in positive characteristic} \hspace{5mm}

In this section we outline the key results from \cite{bmr} that will be used.

\begin{definition} Let $\mathcal{B}=G/B$ be the flag variety, and $T^* \mathcal{B}$ be the Springer resolution. We refer the reader to Section 4.1.2 (and Section 3.1) of \cite{bmr} for the definition of the Springer fiber $\mathcal{B}_e$ in this context. Let $\iota: \mathcal{B}_e \rightarrow T^* \mathcal{B}$ denote the natural inclusion. By the superscript $^{(1)}$ below we mean the Frobenius twist; see Sec 1.1 of \cite{bmr} for an exposition. Let $\Lambda$ be the corresponding weight lattice, for each $\lambda \in \Lambda$ we have a line bundle $\mathcal{O}(\lambda)$ on $G/B$. \end{definition}

The main result proven by Bezrukavnikov, Mirkovic and Rumynin in Theorem 5.3.1 of \cite{bmr} is the following. 

\begin{theorem} $$D^b \text{Mod}^{\text{fg}}_{e, \lambda}(U \mathfrak{g}) \simeq D^b \text{Coh}_{\mathcal{B}_e^{(1)}}(T^* \mathcal{B}^{(1)})$$ \end{theorem}

\begin{definition} Given a module $M \in \text{Mod}^{\text{fg}}_{e, \lambda}(U \mathfrak{g})$, let $\mathcal{F}_M$ be the complex of sheaves that $M$ is mapped to under the above equivalence. We use the symbol $\chi$ to denote both the Euler characteristic of a coherent sheaf, and also the induced function on the Grothendieck group.  \end{definition} 

Assume that $G$ is quasi-simple and let $\check{\alpha}_0$ be the highest coroot. In Section 6.2.7 of \cite{bmr}, Bezrukavnikov, Mirkovic and Rumynin prove the following (please refer to \cite{bmr} for more details).  

\begin{proposition} \label{dimension} Suppose $\mu \in \Lambda$ that $\langle \mu+\rho, \check{\alpha}_0 \rangle < p$.
$$\text{dim}(T_{0 \rightarrow \mu} M) = \chi(\text{Fr}_{\mathcal{B}}^* [i_*\mathcal{F}_M] \otimes \mathcal{O}(\mu + p \rho))$$ \end{proposition}

\begin{remark} In the above formula, we are abusing notation slightly: the image of $\mu$ in $\Lambda/p \Lambda$ gives an element in $\mathfrak{h}^*$, and the subscript in $T_{0 \rightarrow \mu}$ refers to this element. Note also that in \cite{bmr}, we have the statement in the slightly more general setting where $\mu$ lies in the closure of the fundamental alcove. We remark that we follow the conventions from Section 6.2.7 and Lemma 6.2.5 of \cite{bmr}, and in particular the choice of splitting bundle. \end{remark}

\subsection{Affine tangles} The combinatorics of affine tangles will be necessary for what follows; we recall the definitions here (see Section 3 of \cite{an} for more details). 

\begin{definition} If $q \equiv r \pmod 2$, a $(q,r)$ affine tangle is an embedding (up to isotopy) of $\frac{q+r}{2}$ arcs and a finite number of circles into the region $\{ (x, y) \in \mathbb{C} \times \mathbb{R} | 1 \leq |x| \leq 2 \}$, such that the end-points of the arcs are $(1,0), (\zeta_q, 0), \cdots , (\zeta_q^{q-1}, 0), (2, 0), (2\zeta_r,0), \cdots, (2 \zeta_r^{r-1},0)$ in some order; here $\zeta_k = e^{\frac{2\pi i}{k}}$. \end{definition}

\begin{figure} 
\centering
\captionsetup{justification=centering}
\includegraphics[scale=1]{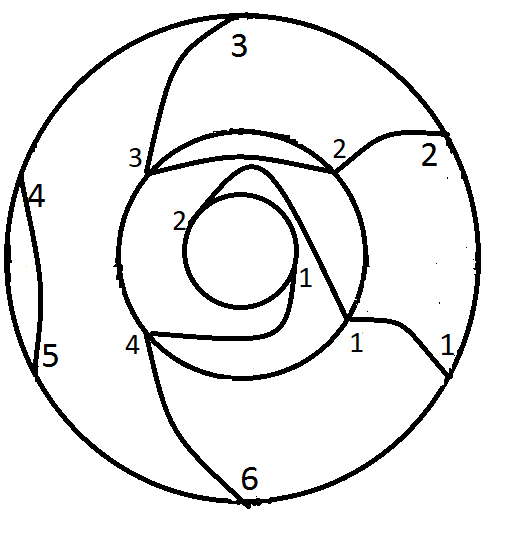}
\caption{}
\label{composition}
\end{figure}

Figure \ref{composition} illustrates the composition of affine tangles.

\begin{definition} If $q \equiv r \pmod 2$, a framed $(q,r)$ affine tangle is an embedding (up to isotopy) of $\frac{q+r}{2}$ ``rectangular arcs" and a finite number of circles into the region $\{ (x, y) \in \mathbb{C} \times \mathbb{R} | 1 \leq |x| \leq 2 \}$. Here a ``rectangular arc" is an injective map from $[0,1] \times [0,1]$ to $\{ (x, y) \in \mathbb{C} \times \mathbb{R} | 1 \leq |x| \leq 2 \}$ where the segments $[0,1] \times 0$ and $[0,1] \times 1$ are mapped to the segments $[x] \times [0,1]$ where $x = \zeta_p^i$ or $x = 2 \zeta_q^j$. 
 \end{definition} 

\begin{remark} Framed tangles are necessary for the validity of the Reidemeister I move in the functorial tangle categorification statement from Proposition \ref{tanglecat}. Below we will often abbreviate ``$(q,r)$-affine tangle'' to ``$(q,r)$-tangle''. In \cite{an}, we used $\; \widehat{} \;$ notation to denote framed tangles; we omit that here. Note that given a $(q,r)$ framed affine tangle $\alpha$, and a $(r,s)$ framed affine tangle $\beta$, we can compose them using scaling and concatenation to obtain a $(q,s)$ affine tangle $\beta \circ \alpha$. This composition is associative up to isotopy. \end{remark}

\begin{definition} Given $1 \leq i \leq n$, define the following affine tangles: \begin{itemize} \item Let $g_n^i$ denote the $(n-2, n)$ tangle with an arc connecting $(2\zeta_n^i,0)$ to $(2\zeta_n^{i+1},0)$.
\item Let $f_n^i$ denote the $(n, n-2)$ tangle with an arc connecting $(\zeta_{n}^i,0)$ and $(\zeta_{n}^{i+1},0)$. 
\item Let $t_n^i(1)$ (respectively, $t^i_n(2)$) denote the $(n,n)$ tangle in which a strand connecting $(\zeta_n^{i},0)$ to $(2\zeta_{n}^{i+1},0)$ passes above (respectively, beneath) a strand connecting $(\zeta_n^{i+1},0)$ to $(2\zeta_n^{i},0)$. 
\item Let $r_n$ denote the $(n,n)$ tangle connecting $(\zeta_n^j, 0)$ to $(2\zeta_n^{j-1}, 0)$ for each $1 \leq j \leq n$ (clockwise rotation of all strands). 
\item Let $s_n^{i}$ denote the $(n,n)$-tangle with a strand connecting $(\zeta_j,0)$ to $(2 \zeta_j,0)$ for each $j$, and a strand connecting $(\zeta_i,0)$ to $(2 \zeta_i,0)$ passing clockwise around the circle, beneath all the other strands.
\end{itemize} 
\end{definition}

We also introduce new generators $w_n^i(1)$ and $w_n^i(2)$, which correspond to positive and negative twists of framing of the $i$-th strand of an $(n, n)$ identity tangle.

The below figures depict some of these elementary tangles. Here the label $i$ in the inner (resp. outer) circle denotes the point $(\zeta^i, 0)$ (resp. $(2 \zeta^i,0)$). 

\includegraphics[scale=0.65]{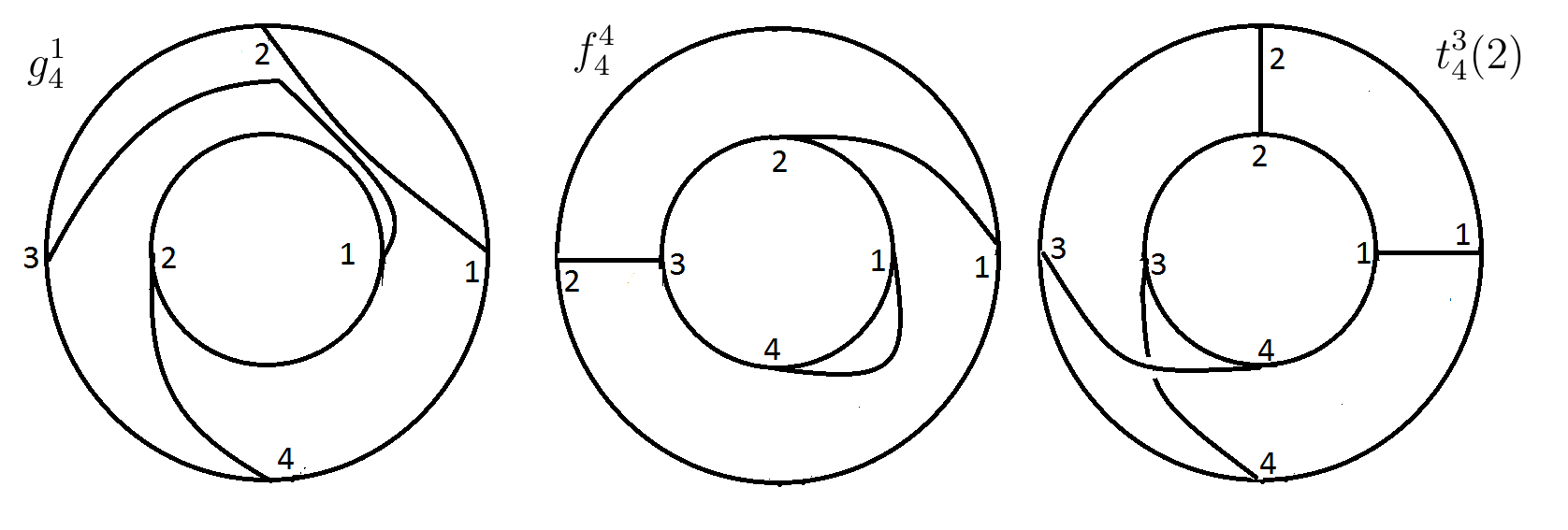}

\includegraphics[scale=0.65]{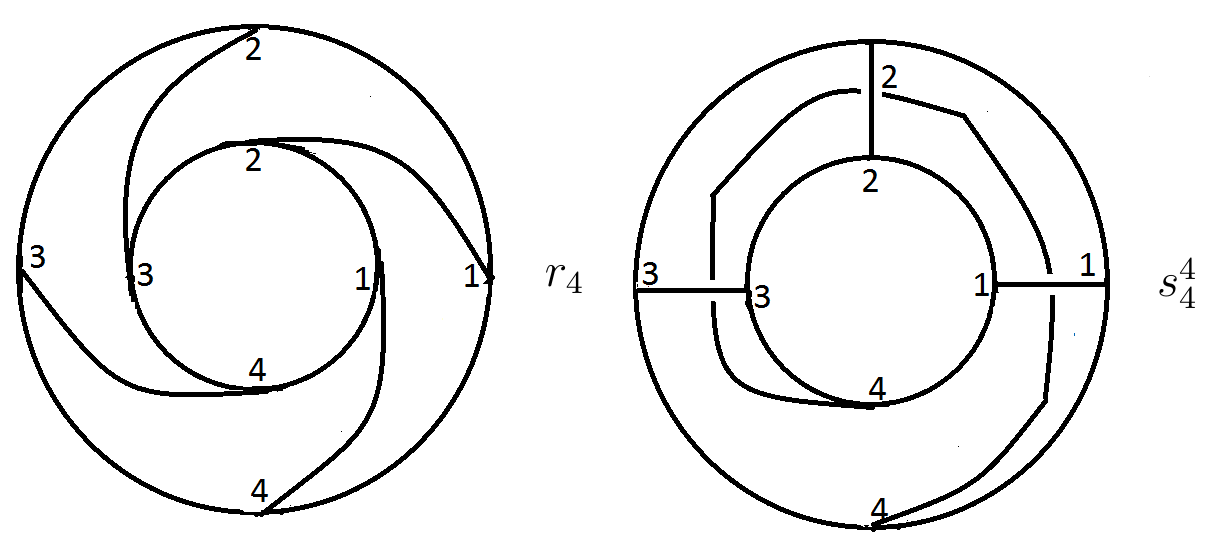}

These are called ``elementary'' tangles because any framed tangle can be expressed as (or rather, is isotopic to) a product of the tangles $g_n^i, , f_n^i, t_n^i(1), t_n^i(2), w_n^i(1), w_n^i(2)$ and $r_n$. Alternatively, this statement is true if we replace $r_n$ by $s_n^n$ (to see this, one uses the fact that $r_n = s_n^n \circ t_n^{n-1}(2) \circ \cdots \circ t_n^1(2)$). 

\begin{definition} A $(q,r)$-tangle is a``linear tangle'' if it can be expressed in terms of the generators $g_n^i, f_n^i, t_n^i(1), t_n^i(2), w_n^i(1), w_n^i(2)$ (i.e. without using $r_n$). \end{definition}

The above definition is motivated that the fact any such tangle can be represented with endpoints on a line segment (as opposed to an annulus); see Section 4 of \cite{ck} for a pictorial depiction.


Below is the list of relations which the generators $g_n^i, f_n^i, t_n^i(1), t_n^i(2), w_n^i(1), w_n^i(2), r_n$ satisfy. Note that $r^{\prime}_n$ is the inverse of $r_n$. 

\begin{enumerate} 
\item\label{AFTanMovesFirst} ${f}_n^i\circ {g}_n^{i+1} = id = {f}_n^{i+1} \circ {g}_n^i$ 
\item\label{AFTanReidemeister1} (Reidemeister 1) ${f}_n^i\circ {t}_n^{i\pm 1}(l)\circ {g}_n^i = {w}_n^i(l)$ 
\item\label{AFTanMovesThird} ${t}_n^i(2)\circ {t}_n^i(1) = id = {t}_n^i(1)\circ {t}_n^i(2)$ 
\item\label{AFTanReidemeister3} ${t}_n^i(l)\circ {t}_n^{i+1}(l)\circ {t}_n^i(l) =  {t}_n^{i+1}(l)\circ {t}_n^{i}(l)\circ {t}_n^{i+1}(l)$ 
\item ${g}_{n+2}^{i+k}\circ {g}_n^i = {g}_{n+2}^i\circ {g}_n^{i+k-2}$ 
\item ${f}_{n}^{i+k-2}\circ {f}_{n+2}^i = {f}_{n}^i\circ {f}_{n+2}^{i+k}$ 
\item ${g}_n^{i+k-2}\circ {f}_n^i = {f}_{n+2}^i\circ {g}_{n+2}^{i+k}, \quad {g}_n^{i}\circ {f}_n^{i+k-2} = {f}_{n+2}^{i+k}\circ {g}_{n+2}^{i}$
\item\label{AFTanPitchfork} ${g}_n^i\circ {t}_{n-2}^{i+k-2}(l) = {t}_n^{i+k}(l)\circ {g}_n^i, \quad {g}_n^{i+k}\circ {t}_{n-2}^{i}(l) = {t}_n^{i}(l)\circ {g}_n^{i+k}$ 
\item ${f}_n^i\circ {t}_{n}^{i+k}(l) = {t}_{n-2}^{i+k-2}(l)\circ {f}_n^i, \quad {f}_n^{i+k}\circ {t}_{n}^{i}(l) = {t}_{n-2}^{i}(l)\circ {f}_n^{i+k}$ 
\item ${t}_n^i(l)\circ {t}_n^{i+k}(m) =  {t}_n^{i+k}(m)\circ {t}_n^i(l)$ \item\label{TanMovesLast}  ${t}_n^i(1)\circ {g}_n^{i+1} = {t}_n^{i+1}(2)\circ {g}_n^i,\quad {t}_n^i(2)\circ {g}_n^{i+1}= {t}_n^{i+1}(1)\circ {g}_n^i$ 
\item ${r}_n\circ {r}_n'  = id = {r}_n' \circ {r}_n$ 
\item ${r}_{n-2}'\circ {f}_n^i \circ {r}_n  = {f}_n^{i+1},\ i=1,\ldots, n-2; \quad
 {f}_n^{n-1} \circ ({r}_n)^2  = {f}_n^1$ 
\item ${r}_n' \circ {g}_n^i \circ {r}_{n-2}  = {g}_n^{i+1}, \ i=1,\ldots, n-2; \quad
 ({r}_n')^{2} \circ {g}_n^{n-1} = {g}_n^1$ 
\item  ${r}_n' \circ {t}_n^i(l)\circ {r}_n  = {t}_n^{i+1}(l); \quad
 ({r}_n')^{2}\circ {t}_n^{n-1}(l) \circ ({r}_n)^2  = {t}_n^1(l)$ 

We have the following additional relations for twists:

\item\label{AFTanMovesFirstTwist} ${w}_n^i(1) \circ {w}_n^i(2) = id, \quad {w}_n^i(l) \circ {w}_n^j(k) = {w}_n^j(k) \circ {w}_n^i(l),\ i\ne j$ 
\item $\hat{w}_n^i(k) \circ {g}_n^i = {w}_n^{i+1}(k) \circ {g}_n^i, \quad%
 {w}_n^i(k) \circ {g}_n^j =  {g}_n^j \circ {w}_n^{i+1\pm 1}(k),\ i\ne j, j+1$
\item ${f}_n^i \circ {w}_n^i(k) = {f}_n^i \circ {w}_n^{i+1}(k), \quad%
{w}_n^i(k) \circ {f}_n^j =  {f}_n^j \circ {w}_n^{i-1\pm 1}(k),\ i\ne j, j+1$
\item ${w}_n^i(k) \circ {t}_n^i(l) = {w}_n^{i+1}(k) \circ {t}_n^i(l), \quad%
 {w}_n^i(k) \circ {t}_n^j(l) =  {t}_n^j(l) \circ {w}_n^{i}(k),\ i\ne j, j+1$ 
\item \label{FTanMovesLast} ${t}_n^i \circ {w}_n^i(k) = {f}_n^i \circ {w}_n^{i+1}(k), \quad%
{w}_n^i(k) \circ {f}_n^j =  {t}_n^j \circ {w}_n^{i}(k),\ i\ne j, j+1$
\item\label{AFTanMovesLast} ${w}_n^i(k) \circ {r}_n = {r}_n \circ {w}_n^{i-1}(k), \quad%
 {w}_n^i(k) \circ {r}'_n = {r}'_n \circ {w}_n^{i+1}(k)$ \end{enumerate}

\subsection{Exotic t-structures for two-block nilpotents} \label{[AN]}

In this section we summarize the results that we will need from \cite{an}. Let $\mathfrak{g}=\mathfrak{sl}_{m+2n}$, and $\mathcal{B}$ be the flag variety. Let $e \in \mathfrak{g}^*$ be a nilpotent linear functional with Jordan type $(m+n, n)$, and $\lambda \in \mathfrak{h}^* // W$ be a regular Harish-Chandra character. 

Let $S_{m,n} \subset \mathfrak{g}$ be the corresponding Mirkovic-Vybornov transverse slice (see \cite{mv} for a definition, and see below for an explicit matrix realization). Let $U_{m+n,n} \subset \tilde{\mathfrak{g}} = T^* \mathcal{B}$ be its preimage under the Grothendieck-Springer resolution map.  

\begin{align*} S_{m,n} &= \{ \left( \begin{array}{ccccccccc}
 & 1 &  & & &  &  &  &  \\
 &  & & \ddots &  &  &  &  &  \\
 &  &  & & 1 &  &  &  &  \\
a_1 & a_2 & & \cdots & a_{m+n} & b_1 & b_2 & \cdots  & b_{n} \\
 &  &  &  & &  & 1 &  &  \\
 &  &  &  & &  &  & \ddots &  \\
 &  &  &  & &  &  &  & 1 \\
c_1 & \cdots & c_n \; 0 & \cdots & 0 & d_1 & d_2 & \cdots & d_n
\end{array} \right) \} \end{align*}

\begin{definition} Let $\langle \cdot, \cdot \rangle$ denote the Killing form. 
\begin{align*} S^*_{m,n} &= \{ x \in \mathfrak{g} \; | \; \langle x, z - e \rangle = 0 \; \forall \; z \in S_{m,n} \} \\ \textbf{U}_{\hat{e}} \mathfrak{g} &= \textbf{U} \mathfrak{g} / \langle x^p - x^{[p]} - \langle e,x \rangle^p \; | \; x \in S^*_{m,n} \rangle \end{align*}

Denote by $\text{Mod}_{\hat{e}, \lambda}^{fg}(\textbf{U}\mathfrak{g})$ the subcategory of $\text{Mod}_{e, \lambda}^{fg}(\textbf{U}\mathfrak{g})$ consisting of those modules which are defined over the quotient $\textbf{U}_{\hat{e}} \mathfrak{g}$. \end{definition}

[BMR] localization theory gives the equivalence in the below Proposition; see Section 5.1 of \cite{an} for more details. Under this equivalence, the tautological t-structure on the LHS is mapped to the exotic t-structure on the RHS. We refer to the images of the irreducible objects in the LHS as simple exotic sheaves.

\begin{proposition} \label{cmn} $$\mathcal{C}_{m,n} := D^b \text{Mod}^{\text{fg}}_{\hat{e}, \lambda}(U\mathfrak{g}) \simeq D^b \text{Coh}_{\mathcal{B}_e}(U_e)$$ \end{proposition}

\begin{proposition} \label{tanglefunctor} Let $\alpha$ be an $(m+2p, m+2q)$-tangle. Then we have a functor $\Psi(\alpha): \mathcal{C}_{m,p} \rightarrow \mathcal{C}_{m,q}$; this collection of functors satisfy the functorial tangle relations $\Psi(\alpha) \circ \Psi(\beta) \simeq \Psi(\alpha \circ \beta)$. \end{proposition}

These are defined by Fourier-Mukai transforms on the categories of coherent sheaves; see Section $4$ of \cite{an} for more details of the proof. These functors can be used to describe the simple exotic sheaves as follows.   
\begin{definition} Each $\alpha \in \text{Cross}(m,n)$ can be viewed as a framed tangle, with the trivial framing. Thus using Proposition \ref{tanglefunctor} we have a functor $\Psi(\alpha): \mathcal{C}_{m,0} \rightarrow \mathcal{C}_{m,n}$. Since $\mathcal{B}_e$ is a point when $n=0$, it follows that $\mathcal{C}_{m,0} \simeq D^b(\text{Vect})$.Given $\alpha \in \text{Cross}(m,n)$, let $\Psi_{\alpha} = \Psi(\alpha)(\textbf{k})$. \end{definition}
\begin{theorem} The objects $\{ \Psi_{\alpha} \; | \; \alpha \in \text{Cross}(m,n) \}$ are the simple exotic sheaves. \end{theorem} 
See Section 5.3 of \cite{an} for a proof of this theorem. Taking $\lambda=0$ in the equivalence in Definition \ref{cmn}, let $M_{\alpha} \in \text{Mod}^{\text{fg}}_{e, 0}(U \mathfrak{g})$ be the irreducible object corresponding to the complex of coherent sheaves $\Psi_{\alpha}$. It is not immediately clear how this description can be matched with that from Proposition 2.7; see however conjecture preceding Remark 5.16 below. 

\subsection{The $\mathfrak{g}=\mathfrak{sl}_2$ example.} Consider the case where $\mathfrak{g}=\mathfrak{sl}_2$, $e=0$ and $\lambda = 0$. The dimensions of the two simple objects in $\text{Mod}^{fg}_{0, \lambda}(U \mathfrak{g})$ are $1$ and $p-1$ (see Section \ref{modularrep} for more details). Under [BMR] localization, they correspond to complexes of coherent sheaves lying in $D^b \text{Coh}(\mathbb{P}^1)$, known as ``simple exotic sheaves''; they are $\mathcal{O}(-1)$ and $\mathcal{O}(-2)[1]$ (see Section 2.5.6 of \cite{bm} for more details, but the normalizations there are different). The formula $\text{dim}(M_{\alpha}) = \chi(\text{Fr}^* \Psi_{\alpha} \otimes \mathcal{O}(p \rho))$ gives us the correct dimensions: \begin{align*} \chi(\text{Fr}^* \mathcal{O}(-1) \otimes \mathcal{O}(p \rho)) &= \chi(\mathcal{O}(-p) \otimes \mathcal{O}(p)) \\ &= \chi(\mathcal{O}) = 1 \\ \chi(\text{Fr}^* \mathcal{O}(-2)[1] \otimes \mathcal{O}(p \rho)) &= \chi(\mathcal{O}(-2p)[1] \otimes \mathcal{O}(p)) \\ &= \chi(\mathcal{O}(-p)[1]) = - \chi(\mathcal{O}(-p)) \\ &= p-1 \end{align*} 

\section{Coherent sheaves on Cautis-Kamnitzer varieties} \label{cohsheaves}

In this section we collect some results about the category of coherent sheaves on the Cautis-Kamnitzer varieties $Y_k$. Since these have not appeared in the literature we give a detailed exposition; but they are all known to experts, and none are new (except for the first part of Lemma \ref{computation}).   

\subsection{Cautis-Kamnitzer's tangle categorification} We will need Cautis-Kamnitzer's framework for much of what follows. The varieties $Y_k$ are constructed as follows (see Section 2.1 of \cite{ck} for more details):

\begin{definition} Consider a $2k$-dimensional vector space $V_k$ with basis $e_1, \ldots , e_{k}, f_1, \ldots ,f_k$ and a nilpotent $z$ such that $z e_{i+1} = e_{i}$, $z f_{i+1} = f_{i}$ for $i \geq 1$, and $ze_1 = zf_1=0$. Let $L_0=0$. Then:  \begin{align*} Y_{k} = \{ (L_1 \subset \cdots \subset L_{k} \subset V_{k}) \; | \; \text{ dim}(L_i) = i,\ zL_{i} \subset L_{i-1} \} \end{align*} 
Denote by $K^0(\text{Coh}(Y_{k}))$ the (complexified) Grothendieck group of $\text{Coh}(Y_k)$, the category of coherent sheaves on $Y_k$.\end{definition}

We have a natural forgetful map $p_k: Y_k \rightarrow Y_{k-1}$ which is a $\mathbb{P}^1$-bundle. To see this, suppose that we have $(L_1,..., L_{k-1}) \in Y_{k-1}$ and are considering possible choices of $L_k$. It
is easy to see that we must have $L_{k-1} \subset L_k \subset z^{-1}(L_k)$. Since $z^{-1}(L_{k-1})/L_{k-1}$ is always two dimensional,
this fibre is a ${\mathbb P}^1$. It follows that $Y_k$ is an iterated $\mathbb{P}^1$-bundle of dimension $k$, since $Y_1 \simeq \mathbb{P}^1$.

\begin{definition} Let $\mathcal{V}_s$ be the tautological vector bundle on $Y_k$ corresponding to $L_s$. For each $1 \leq s \leq k$, we have a line bundle $\Lambda_s = \mathcal{V}_s/\mathcal{V}_{s-1}$ on this space (note that $\Lambda_k = \mathcal{O}(-1)_{Y_k / Y_{k-1}}$, i.e. the relative $\mathcal{O}(-1)$ on the $\mathbb{P}^1$-fibration $Y_k \rightarrow Y_{k-1}$). Given $i_s \in \{ 0,1 \}$ for each $1 \leq s \leq k$, define also the following line bundle: $$\Lambda_{i_1, i_2, \cdots, i_k} := \bigotimes_{1 \leq s \leq k} \Lambda_{s}^{i_s}$$ \end{definition} 

Recall from Section $6.2$ of \cite{ck} that:
\begin{proposition} A basis of $K^0(\text{Coh}(Y_{k}))$ is given by the classes of the line bundles $[\Lambda_{i_1, i_2, \cdots, i_k}]$, where $i_s \in \{ 0, 1 \}$ for each $1 \leq s \leq k$. \end{proposition}

\begin{definition} \label{def:K^0} Consider the vector space $\mathbb{C}^2$, with a fixed basis consisting of the vectors $v_0$ and $v_1$. Identify $K^0(\text{Coh}(Y_{k}))$ with $(\mathbb{C}^2)^{\otimes k}$, so that $$ [\Lambda_{i_1, \cdots, i_{k}}] = (v_1 + \delta_{i_1, 0} v_0) \otimes \cdots \otimes (v_1 + \delta_{i_{k}, 0} v_0) $$ \end{definition} \begin{definition} Denote the monomial basis for $(\mathbb{C}^2)^{\otimes k}$ as follows. Let $I \subset \{ 1, \cdots, k \}$ be a subset, and for $1 \leq a \leq n$ let $I(a) = 1$ if $a \in I$, and $I(a) = 0$ if $a \notin I$. Define: $$v_I = v_{I(1)} \otimes v_{I(2)} \otimes \cdots \otimes v_{I(k)} $$ Let $T_{i_1, \cdots, i_k}$ be the operator on $K^0(\text{Coh}(Y_k))$ induced by tensoring with the line bundle $\Lambda_{i_1, \cdots, i_k}$. \end{definition} 

\begin{remark} While the results are stated and proven in \cite{ck} over the ground field $\mathbb{C}$, they hold more generally over fields of positive characteristic. Another key difference is that in \cite{ck} the results are stated in the equivariant setting; the same proofs can be used to obtain analogous results in the non-equivariant case, and are in fact easier to check here. The above identification of the Grothendieck group is slightly different from that used by Cautis and Kamnitzer in \cite{ck}. Since the class of a point in $K^0(\mathbb{P}^1)$ is $v_0$ with these conventions, it is better suited to our purposes (in particular, both Lemmas \ref{computation} and \ref{base} are cleaner with this choice of notation). \end{remark}  

\begin{definition} Let $\mathcal{D}_k = D^b\text{Coh}(Y_k)$ be the derived category of coherent sheaves on $Y_k$. \end{definition}

The below proposition is a strengthening of the results in Section 6 of \cite{ck}, the difference being that we use affine tangles (and not linear tangles). See that section for the definition of Reshetikhin-Turaev invariant $\psi(\alpha): (\mathbb{C}^2)^{\otimes q} \rightarrow (\mathbb{C}^2)^{\otimes r}$; here we have specialized the parameter of the quantum group $U_v(\mathfrak{sl}_2)$ to $v=1$, so that we are working with representations of $\mathfrak{sl}_2$ (and not $U_v(\mathfrak{sl}_2)$). 

\begin{definition} For each linear $(q,r)$-tangle $\alpha$, we define the map $\psi(\alpha): (\mathbb{C}^2)^{\otimes q} \rightarrow (\mathbb{C}^2)^{\otimes r}$ as follows. It suffices to define the map on the elementary linear tangles. Below $k \in \{ 1, 2 \}$ and $i,j \in \{0, 1 \}$.
\begin{align*}\psi(f_2^1): & (\mathbb{C}^2)^{\otimes 2} \rightarrow \mathbb{C}; \qquad \psi(f_2^1)(v_i \otimes v_i) = 0, \psi(f_2^1)(v_0 \otimes v_1) = 1, \psi(f_2^1)(v_1 \otimes v_0) = -1 \\ \psi(g_2^1): & \mathbb{C} \rightarrow (\mathbb{C}^2)^{\otimes 2}; \qquad \psi(g_2^1)(1) = v_1 \otimes v_0 - v_0 \otimes v_1 \\ \psi(t_2^1(k)): & (\mathbb{C}^2)^{\otimes 2} \rightarrow (\mathbb{C}^2)^{\otimes 2}; \qquad \psi(t_2^1(k)) (v_i \otimes v_j) = v_j \otimes v_i \\ \end{align*} 
Now we define $\psi(g_n^i)$ as follows, and similarly for $\psi(f_n^i)$ and $\psi(t_n^i(k))$. 
\begin{align*} \psi(g_n^i): & (\mathbb{C}^2)^{\otimes n-2} \rightarrow (\mathbb{C}^2)^{\otimes n};  \psi(g_n^i) = \text{id}_{(\mathbb{C}^2)^{\otimes i-1}} \otimes \psi(g_2^1) \otimes  \text{id}_{(\mathbb{C}^2)^{\otimes n-2-i+1}}  \end{align*}
Note that $\psi$ is well-defined due to \cite{rt} and the introduction to \cite{ck}.
\end{definition} 

\begin{proposition} \label{tanglecat} For each $(q,r)$-tangle $\alpha$, there exist a functor $\Psi(\alpha): \mathcal{D}_q \rightarrow \mathcal{D}_r$, and these are compatible under composition, i.e. given an $(r,s)$-angle $\beta$, $\Psi(\alpha \circ \beta) \simeq \Psi(\alpha) \circ \Psi(\beta)$. When $\alpha$ is a linear tangle, the image of $\Psi(\alpha)$ in the Grothendieck group is the Reshetikhin-Turaev invariant $\psi(\alpha)$. \end{proposition} \begin{proof} When $\alpha$ is a linear tangle, the functors are certain Fourier-Mukai kernels which are described in Section 4.2 of \cite{ck}; the functorial relations are proven in Section 5. Theorem 6.5 of \cite{ck} proves that, on the level of Grothendieck groups, $\Psi(\alpha)$ corresponds to $\psi(\alpha)$ (but one needs to take into account that our identification of $K^0(Y_k)$ with $(\mathbb{C}^2)^{\otimes k}$ differs from theirs). It remains to construct the functors $\Psi(\alpha)$ when $\alpha$ is an affine tangle; this follows from the techniques used in Section 4.5 of \cite{an}. It suffices to construct a functor $\Psi(s_{k}^k)$, and check that it satisfies the necessary functorial tangle relations (13) and (14) on page 7. Define $\Psi(s_{k}^k)$ to be the functor of tensoring with the line bundle $\Lambda_k^{-1}$; to check the functorial tangle relations, the same argument used in Theorem 4.27 and Proposition 4.25 of \cite{an} works verbatim. \end{proof}

\begin{remark} For the above proposition, ribbon framings are necessary for the validity of the Reidemeister I move: $f_n^i \circ t_n^{i \pm 1}(1) \circ g_n^i = w_n^i(1)$. On the functorial level, $w_n^i(1)$ corresponds to a homological shift. \end{remark} 

\subsection{Embedding of Slodowy slices} \label{embedding}

\begin{definition} Let $W_{m,n} \subset V_{m+2n}$ denote the vector subspace with basis $e_1, \ldots, e_{m+n}, f_1, \cdots, f_n$, so that $z|_{W_{m,n}}$ has Jordan type $(m+n,n)$. Identify $\mathfrak{sl}_{m+2n}$ with $\mathfrak{sl}(W_{m,n})$. Let $P:V_{m+2n} \rightarrow W_{m,n}$ denote the projection defined by $Pe_{i} = e_{i}$ if $i \leq m+n$, $Pe_{i}=0$ if $i > m+n$; $Pf_i = f_i$ if $i \leq n$, $Pf_i = f_i$ if $i > n$. Denote by $\mathcal{B}_{m+n,n}$ the Springer fiber for $z \in \mathfrak{sl}(W_{m,n})$. $$ \mathcal{B}_{m+n,n} = \{ (0 \subset U_1 \subset \cdots \subset U_{m+2n-1} \subset U_{m+2n} = W_{m,n}) \; | \; z U_i \subseteq U_{i-1} \} $$ \end{definition}

\begin{proposition} \label{embed} There exists a locally closed embedding $i: U_{m+n,n} \rightarrow Y_{m+2n}$. The embedding $i$ is the identity on the Springer fiber $\mathcal{B}_{m+n,n}$.
\end{proposition}

In this Section we prove the above Proposition, mirroring the approach used to prove Proposition 2.4 in \cite{ck} (this fact is known in greater generality using results of Mirkovic-Vybornov, \cite{mv}). The proof will not be used elsewhere in the paper, and the reader may wish to skip the rest of this subsection. In the below definition, note that $P(L_{m+2n})=W_{m,n}$ is an open condition since it is equivalent to the coefficients of a basis of $L_{m+2n}$ being nonzero, hence $Q_{m,n}$ is an open subvariety. 

\begin{definition} Define the open subvariety $Q_{m,n} \subseteq Y_{m+2n}$ as: $$Q_{m,n} =\{ (L_1 \subset \cdots \subset L_{m+2n}) \in Y_{m+2n} \; | \; P(L_{m+2n}) = W_{m,n} \}$$ \end{definition}

\begin{definition} Define the affine variety $\overline{S}'_{m,n}$ to consist of all matrices of the following form. \begin{align*} \overline{S}_{m,n}' &= \{ \left( \begin{array}{cccccccc}
 & 1 &  &  &  &  &  &  \\
 &  & \ddots &  &  &  &  &  \\
 &  &  & 1 &  &  &  &  \\
a_1 & a_2 & \cdots & a_{m+n} & b_1 & b_2 & \cdots  & b_{n} \\
 &  &  &  &  & 1 &  &  \\
 &  &  &  &  &  & \ddots &  \\
 &  &  &  &  &  &  & 1 \\
c_1 & c_2 & \cdots & c_{m+n} & d_1 & d_2 & \cdots & d_n
\end{array} \right) \} \end{align*}
 \end{definition}

\begin{definition} Define $S'_{m,n}$ to be the subvariety of $\overline{S}_{m,n}'$ consisting of all nilpotent matrices. \end{definition}

Note that we have $U_{m+n,n} \subset S_{m,n}' \times_{\mathfrak{sl}_{m+2n}} T^* \mathcal{B}$.

\begin{proof}[Proof of Proposition \ref{embed}] The existence of the map $i: U_{m+n,n} \rightarrow Y_{m+2n}$ (and the description of $i|_{\mathcal{B}_{m+n,n}}$) follows from the below Lemma, via the composition: $$ U_{m+n,n} \hookrightarrow S_{m,n}' \times_{\mathfrak{sl}_{m+2n}} T^* \mathcal{B} \simeq Q_{m,n} \hookrightarrow Y_{m+2n} $$ \end{proof} 

\begin{lemma} \label{sub} 
There is an isomorphism $Q_{m,n} \simeq S_{m,n}' \times_{\mathfrak{sl}_{m+2n}} T^* \mathcal{B}$. 
\end{lemma}
\begin{proof}[Proof of Lemma \ref{sub}]
Given $(L_1 \subset \cdots \subset L_{m+2n}) \in Q_{m,n}$, since $P: L_{m+2n} \rightarrow W_{m,n}$ is an isomorphism, we have a nilpotent endomorphism $x = PzP^{-1} \in \text{End}(W_{m,n})$. If $P^{-1}e_i = e_{i} + v'$, where $v'$ lies in the span of $e_{m+n+1}, \cdots, e_{m+2n}, f_{n+1}, \cdots, f_{m+2n}$, then $zP^{-1}e_i=e_{i-1}+v''$ where $v''$ is in the span of $e_{m+n}, \cdots, e_{m+2n-1}, f_n, \cdots, f_{m+2n-1}$. Hence $PzP^{-1}e_i = x e_{i} \in e_{i-1} + \text{span}(e_{m+n}, f_n)$, and similarly $x f_{i} \in f_{i-1} + \text{span}(e_{m+n}, f_n)$; so $x \in S_{m,n}'$. Thus we have a map $\alpha: Q_{m,n} \rightarrow S_{m,n}' \times_{\mathfrak{sl}_{m+2n}} T^* \mathcal{B}$ given by $\alpha(L_1, \cdots, L_{m+2n}) = (PzP^{-1}, (P(L_1), P(L_2), \cdots, P(L_{m+2n})))$.

For the converse direction, from the below Lemma \ref{lin} we know that given a nilpotent $x \in S_{m,n}'$ there exists a unique $z$-stable subspace $L_{m+2n} \subset V_{m+2n}$ such that $P L_{m+2n} = W_{m,n}$ and $PzP^{-1} = x$; call this subspace $L_{m+2n} = \Theta(x)$. We have an isomorphism $P: \Theta(x) \simeq W_{m,n}$. Thus given an element $((0 \subset V_1 \subset \cdots \subset V_{m+2n}), x) \in S_{m,n}' \times_{\mathfrak{sl}_{m+2n}} T^* \mathcal{B}$, let $\beta(x)=(0 \subset P^{-1}V_1 \subset P^{-1}V_2 \subset \cdots \subset \Theta(x))$. It is clear that $\alpha$ and $\beta$ are inverse to one another. Note also that the variety $S_{m,n}' \times_{\mathfrak{sl}_{m+2n}} T^* \mathcal{B}$ is smooth, since it is the resolution of $S_{m,n}'$. Any bijective morphism of algebraic varieties with smooth image is necessarily an isomorphism, completing the proof. 
\end{proof}

\begin{lemma} \label{lin} Given a nilpotent $x \in S_{m,n}'$, there exists a unique subspace $L_{m+2n} \subset V_{m+2n}$, with $P L_{m+2n} = W_{m,n}$, such that $z L_{m+2n} \subset L_{m+2n}$ and $PzP^{-1}=x$. \end{lemma}
\begin{proof}[Proof of Lemma \ref{lin}] Since $P L_{m+2n} = W_{m,n}$, to specify the subspace $L_{m+2n}$ it suffices to specify \begin{align*} \tilde{e}_i:=P^{-1}(e_i) = e_{i} + \sum_{1 \leq k \leq n} a_i^{(k)} e_{m+n+k} + \sum_{1 \leq l \leq m+n} c_i^{(l)}f_{n+l} \\ \tilde{f}_j:=P^{-1}(f_j) = f_j + \sum_{1 \leq k \leq n} b_j^{(k)} e_{m+n+k} + \sum_{1 \leq l \leq m+n} d_j^{(l)} f_{n+l} \end{align*} Suppose for $1 \leq i \leq m+n, 1 \leq j \leq n$, $x e_{i} = e_{i-1} + a_i e_{m+n} + c_i f_{n}, xf_j = f_{j-1} + b_j e_{m+n} + d_j f_{n}$; then the identity $PzP^{-1}=x$ is equivalent to $a_i^{(1)}=a_i, c_i^{(1)} = c_i, b_j^{(1)}=b_j$ and $d_j^{(1)}=d_j$. Assuming that $PzP^{-1}=x$, the statement $z L_{m+2n} \subset L_{m+2n}$ is equivalent to saying that $z \tilde{e}_i, z \tilde{f}_j \in L_{m+2n}$: \begin{align*} z \tilde{e}_i = \tilde{e}_{i-1} + a_i \tilde{e}_{m+n} + c_i \tilde{f}_n \\ z \tilde{f}_j = \tilde{f}_{j-1} + b_j \tilde{e}_{m+n} + d_j \tilde{f}_n \end{align*} 
Expanding the above two equations, we obtain: 
\footnotesize
\begin{multline*}
 e_{i-1} + \sum_{1 \leq k \leq n} a_i^{(k)} e_{m+n+k-1} + \sum_{1 \leq l \leq m+n} c_i^{(l)}f_{n+l-1} = e_{i-1} + \sum_{1 \leq k \leq n} a_{i-1}^{(k)} e_{m+n+k} + \sum_{1 \leq l \leq m+n} c_{i-1}^{(l)}f_{n+l} +\\
 + a_i\left(e_{m+n} + \sum_{1 \leq k \leq n} a_{m+n}^{(k)} e_{m+n+k} + \sum_{1 \leq l \leq m+n} c_{m+n}^{(l)}f_{n+l}\right)
 + c_i\left(f_n + \sum_{1 \leq k \leq n} b_n^{(k)} e_{m+n+k} + \sum_{1 \leq l \leq m+n} d_n^{(l)} f_{n+l}\right);
\end{multline*}
\begin{multline*}
f_{j-1} + \sum_{1 \leq k \leq n} b_j^{(k)} e_{m+n+k-1} + \sum_{1 \leq l \leq m+n} d_j^{(l)} f_{n+l-1} = f_{j-1} + \sum_{1 \leq k \leq n} b_{j-1}^{(k)} e_{m+n+k} + \sum_{1  \leq l \leq m+n} d_{j-1}^{(l)} f_{n+l}+ \\
 + b_j\left(e_{m+n} + \sum_{1 \leq k \leq n} a_{m+n}^{(k)} e_{m+n+k} + \sum_{1 \leq l \leq m+n} c_{m+n}^{(l)}f_{n+l}\right)
 + d_j\left(f_n + \sum_{1 \leq k \leq n} b_n^{(k)} e_{m+n+k} + \sum_{1 \leq l \leq m+n} d_n^{(l)} f_{n+l}\right).
\end{multline*} 
\normalsize
Extracting coefficients of $e_{m+n+k}$ and $f_{n+l}$ in the above two equations gives: \begin{align*} a_i^{(k+1)} &= a_{i-1}^{(k)} + a_i a_{m+n}^{(k)} + c_i b_n^{(k)}, \qquad b_j^{(k+1)} = b_{j-1}^{(k)} + b_j a_{m+n}^{(k)} + d_j b_n^{(k)} \\ c_i^{(l+1)} &= c_i^{(l)} + a_i c_{m+n}^{(l)} + c_i d_n^{(l)}, \qquad d_j^{(l+1)} = d_{j-1}^{(l)} + b_j c_{m+n}^{(l)} + d_j d_n^{(l)} \end{align*} Consider the matrix coefficients $(x^k)_{p,q}$ for $1 \leq p,q \leq m+2n$. It follows by induction that we have $a_i^{(k)}= (x^k)_{m+n,i}, b_j^{(k)}= (x^k)_{m+n,m+n+j}, c_i^{(l)}= (x^l)_{m+2n,i}, d_j^{(l)}=(x^l)_{m+2n,m+n+j}$. Indeed, the case where $k=l=1$ is clear; and the induction step follows from expanding the equation $(x^{r+1})_{uv} = \sum_{1 \leq w \leq m+2n} (x^r)_{uw} (x)_{wv}$ for $u=m+n$ and $u=m+2n$.

Using the above recursive definition of $a_{i}^{(k)}$, $b_{j}^{(k)}$, $c_{i}^{(l)}$, and $d_{j}^{(l)}$, it remains to prove that $a_{i}^{(n+1)}=b_{j}^{(n+1)}=0$ and $c_{i}^{(m+n+1)}=d_{j}^{(m+n+1)}=0$. Thus we must show that $(x^{n+1})_{m+n,p}=(x^{m+n+1})_{m+2n,p}=0$ given $1 \leq p \leq m+2n$. Using the equation $(x^{r+1})_{uv} = \sum_{1 \leq w \leq m+2n} (x)_{uw} (x^r)_{wv}$, we compute that: \begin{align*} (x^{n+1})_{m+n,p}&=(x^{n+2})_{m+n-1,p}= \cdots = (x^{m+2n})_{1,p}=0 \\ (x^{m+n+1})_{m+2n,p}&=(x^{m+n+2})_{m+2n-1,p}= \cdots = (x^{m+2n})_{m+n+1,p}=0 \end{align*} This completes the proof of the existence and uniqueness of a $z$-stable subspace $L_{m+2n} \subset V_{m,n}$ with $P L_{m+2n} = W_{m,n}$ and $PzP^{-1} = x$. \end{proof}

\subsection{Grothendieck group of coherent sheaves}
For our computations, we will need the following lemma which tells us what the effect of tensoring by $\Lambda_i^{-1}$ on $K^0(Y_k)$ is, when $i \leq k$. It more or less follows from the argument used to prove Theorem $6.2$ of \cite{ck} (see equation (20) in particular). 
\begin{lemma} \label{inverse} We have the following identity. $$[\Lambda_{i}^{-1}] = 2 [\mathcal{O}] - [\Lambda_i]$$ \end{lemma} \begin{proof} Consider the two exact sequences below. Note here that $z^{-1} \mathcal{V}_{i-1}$ denotes the vector bundle whose fiber at a point is the pre-image of the corresponding fiber of $\mathcal{V}_{i-1}$ at that point, under the map $z$.  \begin{align*} 0 \rightarrow \text{ker}(z) \rightarrow z^{-1} & \mathcal{V}_{i-1} \rightarrow \mathcal{V}_{i-1} \rightarrow 0 \\ 0 \rightarrow \mathcal{V}_i / \mathcal{V}_{i-1} \rightarrow z^{-1} \mathcal{V}_{i-1} & / \mathcal{V}_{i-1} \rightarrow z^{-1} \mathcal{V}_{i-1} / \mathcal{V}_i \rightarrow 0 \end{align*} Taking determinants, we obtain the following (noting that $\text{ker}(z) \simeq \mathcal{O}^{\oplus 2}, \text{det}(\text{ker}(z)) \simeq \mathcal{O}$): 
\begin{align*} \text{det}(z^{-1} \mathcal{V}_{i-1}) &\simeq \text{det}(\mathcal{V}_{i-1}) \\ \Lambda_i \otimes z^{-1} \mathcal{V}_{i-1} / \mathcal{V}_i &\simeq \text{det}(z^{-1} \mathcal{V}_{i-1} / \mathcal{V}_{i-1}) \\ &\simeq \text{det}(z^{-1} \mathcal{V}_{i-1}) \otimes {\text{det}(\mathcal{V}_{i-1})}^{-1} \simeq \mathcal{O} \\ z^{-1} \mathcal{V}_{i-1} / \mathcal{V}_i &\simeq \Lambda_i^{-1} \end{align*} Using the first exact triangle, $[z^{-1} \mathcal{V}_{i-1} / \mathcal{V}_{i-1}] = 2 [\mathcal{O}]$; so $[\Lambda_{i}^{-1}] = 2 [\mathcal{O}] - [\Lambda_i]$, as required. \end{proof}

\begin{definition} Given an $(m+2p, m+2q)$-tangle $\alpha$, denote by $\overline{\alpha}: (\mathbb{C}^2)^{\otimes m+2p} \rightarrow (\mathbb{C}^2)^{\otimes m+2q}$ the map induced on the Grothendieck group by the functor $\Psi(\alpha): \mathcal{D}_{m+2p} \rightarrow \mathcal{D}_{m+2q}$. Given an element $\alpha \in \text{Cross}(m,n)$, denote by $\underline{\alpha} \in (\mathbb{C}^2)^{\otimes m+2n}$ the image of the irreducible object $\Psi_{\alpha}$ in the Grothendieck group. Denote by $\theta_{\lambda}: (\mathbb{C}^2)^{\otimes m+2n} \rightarrow (\mathbb{C}^2)^{\otimes m+2n}$ the image of tensoring with the line bundle $\mathcal{O}(\lambda)$ in the Grothendieck group, for $\lambda \in \Lambda$. Let $\lambda = \lambda_1 \epsilon_1 + \cdots + \lambda_n \epsilon_n$, where $\epsilon_i \in \mathfrak{h}^*$ picks out the $i$-th co-ordinate; and let $\theta_i := \theta_{\epsilon_i}$. \end{definition}  

\begin{definition} Define the endomorphism $\theta$ of the two-dimensional vector space $\mathbb{C}v_0 \oplus \mathbb{C}v_1$ by $$ \theta(v_0) = v_0; \qquad \theta(v_1) = v_1 - v_0 $$ \end{definition} 
We have $K^0(Y_1) \simeq \mathbb{C}v_0 \oplus \mathbb{C}v_1$, and $\theta$ is the endomorphism induced on the Grothendieck group by tensoring with the line bundle $\Lambda_1$. We will need the following explicit descriptions of the maps $\overline{r_k}$ and $\theta_{\lambda}$:
\begin{lemma} \label{computation} Let $i_s \in \{0,1\}$ for $1 \leq s \leq k$. Then we have:
\begin{align*}  \overline{r_k}(v_{i_1} \otimes \cdots \otimes v_{i_k}) &= v_{i_2} \otimes v_{i_3} \otimes \cdots \otimes v_{i_k} \otimes \theta^{-1}(v_{i_1}) \\ [\Lambda_1^{i_1} \otimes \Lambda_2^{i_2} \otimes \cdots \otimes \Lambda_k^{i_k}] &= (v_1 - (i_1-1) v_0) \otimes \cdots \otimes (v_1 - (i_k-1) v_0) \end{align*} \end{lemma}

\begin{proof}[Proof of Lemma \ref{computation}] Using Lemma \ref{inverse} and Proposition \ref{tanglecat}, we compute that: \begin{align*} \overline{r_k}(v_{i_1} \otimes \cdots \otimes v_{i_k}) &= \overline{s_k^k} \circ \overline{t_k^{k-1}(1)} \circ \cdots \circ \overline{t_k^1(1)} (v_{i_1} \otimes \cdots \otimes v_{i_k}) \\ &= \overline{s_k^k} (v_{i_2} \otimes \cdots \otimes v_{i_k} \otimes v_{i_1}) \\ &= v_{i_2} \otimes \cdots \otimes v_{i_k} \otimes \theta^{-1}(v_{i_1}) \end{align*} It follows from Lemma \ref{inverse} that $[\Lambda_i^l] = l [\Lambda_i] - (l-1) [\mathcal{O}]$ for all $l \in \mathbb{Z}$. From Definition \ref{def:K^0}, it follows that: \begin{align*} [\Lambda_1^{i_1} \otimes \Lambda_2^{i_2} \otimes \cdots \otimes \Lambda_k^{i_k}] &= \theta^{i_1}(v_0+v_1) \otimes \cdots \otimes \theta^{i_k}(v_0+v_1) \\ &= (v_1 - (i_1-1) v_0) \otimes \cdots \otimes (v_1 - (i_k-1) v_0) \end{align*} \end{proof} 

The Euler characteristic $\chi$ can be viewed as a function on the complexified Grothendieck group. Now we show that $\chi$ simply picks out the coefficient of $v_0 \otimes \cdots \otimes v_0$. $$\chi: K^0(Y_k) \simeq (\mathbb{C}^2)^{\otimes k} \rightarrow \mathbb{C}$$ 
\begin{lemma} \label{Eulerchar} $\chi(v_{i_1} \otimes \cdots \otimes v_{i_k})= \delta_{i_1, 0} \delta_{i_2,0} \cdots \delta_{i_k, 0}$ \end{lemma}
 
\begin{proof}[Proof of Lemma \ref{Eulerchar}] Recall the $\mathbb{P}^1$-fibration $p_k: Y_{k} \rightarrow Y_{k-1}$. Note that: \begin{align*}  \chi(\Lambda_{i_1, i_2, \cdots, i_k}) &= \chi(Rp_{k*} \Lambda_{i_1, i_2, \cdots, i_k}) \\ &= \chi(Rp_{k*}(p_k^* \Lambda_{i_1, \cdots, i_{k-1}} \otimes \Lambda_k^{i_k})) \\ &= \chi(\Lambda_{i_1, \cdots, i_{k-1}} \otimes Rp_{k*}(\Lambda_k^{i_k})) \end{align*} If $i_k=1$, then $\Lambda_k$ is isomorphic to $\mathcal{O}(-1)$ on each of the $\mathbb{P}^1$-fibers of $p_k$, and consequently $Rp_{k*}(\Lambda_k^{i_k})=0$, and $\chi(\Lambda_{i_1, i_2, \cdots, i_k})=0$. Otherwise if $i_k = 0$, then $\chi(\Lambda_{i_1, i_2, \cdots, i_k})=\chi(\Lambda_{i_1, i_2, \cdots, i_{k-1}})$. By induction, we have that $\chi(\Lambda_{i_1, i_2, \cdots, i_k})= \delta_{i_1, 0} \delta_{i_2, 0} \cdots \delta_{i_k, 0}$. Using Definition \ref{def:K^0}: $$ \chi( (v_1 + \delta_{i_1, 0}v_0) \otimes \cdots \otimes (v_1 + \delta_{i_k,0} v_0) ) = \delta_{i_1, 0} \delta_{i_2, 0} \cdots \delta_{i_k, 0} $$ It follows that $\chi(v_{i_1} \otimes \cdots \otimes v_{i_k}) = \delta_{i_1, 0} \cdots \delta_{i_k,0}$ \end{proof}

\section{Dimension formulae} \label{dimensions}

Before proving the dimension formulae we will need some more notation. Recall that $C(\alpha)$ was defined in Definition \ref{Calpha}.

\begin{definition} Let $T(\alpha)$ consist of all $n$-element subsets $S \subset \{ 1, \cdots, m+2n \}$ such that $|S \cap C| = 1$ for each $C \in C(\alpha)$; clearly $|T(\alpha)| = 2^n$. Let $\text{max}(C)$ be the maximal element in $C$, and for each $S \in T(\alpha)$ define: $$\text{sgn}(S) = (-1)^{c(S)} \text{ where } c(S) = \# \{ C : C \in C(\alpha), |C| = 2, \; S \cap C= \text{max}(C) \}.$$ \end{definition}

\begin{example} Let $\alpha \in \text{Cross}(2,2)$ be as in Example \ref{cross22}. Then we have:
\begin{align*}  T(\alpha) &= \{ (1,3), (1,6), (2,3), (2,6) \} \\ \text{sgn}(1,3)&=\text{sgn}(2,6)=1, \text{sgn}(1,6)=\text{sgn}(2,3)=-1 \end{align*} \end{example}

\begin{lemma} \label{base} Let $\alpha \in \text{Cross}(m,0)$ (ie. $\alpha$ is the ``identity'' $(m,m)$-tangle). Then: $$\underline{\alpha} = v_0 \otimes \cdots \otimes v_0 \in V^{\otimes m} $$ \end{lemma} \begin{proof} The Springer fiber $\mathcal{B}_{m, 0}$ is a point; we will prove the stronger statement that given any point $y_m \in Y_m$, its class $[\mathcal{O}_{y_m}] \in K^0(Y_m)$ is $v_0 \otimes \cdots \otimes v_0$ by induction. First note that when $m=1$, $Y_1 \simeq \mathbb{P}^1$, and hence $\underline{\alpha}=v_0$ using the following short exact sequence (noting that $\Lambda_1 \simeq \mathcal{O}(-1)$ and $\Lambda_0 \simeq \mathcal{O}$): \begin{align*} 0 \rightarrow \Lambda_{1} \rightarrow &\Lambda_{0} \rightarrow i_* \mathcal{O}_{y_1} \rightarrow 0 \end{align*} 

For the induction step, consider the $\mathbb{P}^1$-bundle $p_m: Y_m \rightarrow Y_{m-1}$. Let $p_m(y_m) = y_{m-1}$, and for every $m$ denote the inclusion of the point by $j_m: \{y_m\} \rightarrow Y_m$. Taking a relative version of the above exact sequence on the fiber gives: \begin{align*} 0 \rightarrow& \; \mathcal{O}_{p_m^{-1}(y_{m-1})} \otimes \mathcal{O}(-1)_{Y_m/Y_{m-1}} \rightarrow \mathcal{O}_{p_m^{-1}(y_{m-1})} \rightarrow j_{m*} \mathcal{O}_{y_m} \rightarrow 0 \end{align*} By the inductive hypothesis, $[\mathcal{O}_{y_{m-1}}]=v_0 \otimes \cdots \otimes v_0$. Since $p_m^* \Lambda_{i_1, \cdots, i_{m-1}} = \Lambda_{i_1, \cdots, i_{m-1}, 0}$, using Definition \ref{def:K^0} we compute that: \begin{align*} [p_m^*]v_{i_1} \otimes \cdots \otimes v_{i_{m-1}} &= v_{i_1} \otimes \cdots \otimes v_{i_{m-1}} \otimes (v_1+v_0) \\ [\mathcal{O}_{p_m^{-1}(y_{m-1})}] &= [p_m^* \mathcal{O}_{y_{m-1}}] = v_0 \otimes \cdots \otimes v_0 \otimes (v_1 + v_0)  \end{align*} Since $\mathcal{O}(-1)_{Y_m/Y_{m-1}} \simeq \Lambda_{0,\cdots,0,1}$, we deduce that $[\mathcal{O}_{p_m^{-1}(y_{m-1})} \otimes \mathcal{O}(-1)_{Y_m/Y_{m-1}}] = v_0 \otimes \cdots \otimes v_0 \otimes v_1$ (as tensoring with $\Lambda_{0,\cdots,0,1}$ corresponds to applying $\theta$ to the last co-ordinate). This now completes the proof: $$[j_{m*} \mathcal{O}_{y_m}] = [\mathcal{O}_{p_m^{-1}(y_{m-1})}] - [\mathcal{O}_{p_m^{-1}(y_{m-1})} \otimes \mathcal{O}(-1)_{Y_m/Y_{m-1}}] = v_0 \otimes \cdots \otimes v_0$$ \end{proof} 

\begin{definition} We say that a matching $\alpha \in \text{Cross}(m,n)$ is ``good'' if $C_1(\alpha) = \emptyset$. \end{definition}

\begin{lemma} \label{class_easy} If the matching $\alpha$ is good, then: \begin{align*} \underline{\alpha} &= \sum_{I \in T(\alpha)} \text{sgn}(I) v_{I} \qquad \qquad \end{align*} \end{lemma}

\begin{proof}[Proof of Lemma \ref{class_easy}] Proof by induction. For the induction step, note that if $\alpha \in \text{Cross}(m,n)$ is a good matching, we can pick $i$ such that $1 < i <m+2n$ and $(i, i+1) \in C(\alpha)$, express $\alpha = g_{m+2n}^i \beta$, for some $\beta \in \text{Cross}(m,n-1)$. Then we can complete the proof using Proposition \ref{tanglecat}: \begin{align*} \underline{\alpha} &= \psi(g_{m+2n}^i) \underline{\beta} \\ &= \sum_{I \in T(\beta)} \text{sgn}(I) \psi(g_{m+2n}^i) v_I  \\ &= \sum_{I \in T(\alpha)} \text{sgn}(I) v_I \end{align*} \end{proof}  

Recall also Lemma 6.2.6 from \cite{bmr}:
\begin{lemma} \label{6.2.6} Given $\mathcal{F} \in D^b(\text{Coh}(\mathcal{B}))$, there exists a polynomial $d_{\mathcal{F}} \in \mathbb{Q}[\Lambda]$, which has the following properties: \begin{align*} \chi(\mathcal{F} \otimes \mathcal{O}(\lambda)) &=d_{\mathcal{F}}(\lambda) \\ d_{\text{Fr}^* \mathcal{F}}(\mu) &= p^{\text{dim}(\mathcal{B})} d_{\mathcal{F}} (\frac{\mu + (1-p) \rho}{p}) \qquad \qquad \blacksquare \end{align*} \end{lemma} 

 
\begin{definition} Given $\alpha \in \text{Cross}(m,n)$, define the polynomial $d_{\alpha}$ via: $$ d_{\alpha}(\tilde{\lambda}) = d_{\Psi_{\alpha}}(\lambda) = \chi(\Psi_{\alpha} \otimes \mathcal{O}(\lambda)) $$ \end{definition}

\begin{lemma} \label{Eulertensoring} \begin{align*} 
d_{\alpha}(\tilde{\lambda}) &=  \prod_{(i,j) \in C_0(\alpha)} (\tilde{\lambda}_i - \tilde{\lambda}_j) \prod_{(i,j) \in C_1(\alpha)} (1 - \tilde{\lambda}_i + \tilde{\lambda}_j) \end{align*} \end{lemma} \begin{proof} 
Let $\alpha = r_{m+2n}^k \beta$ where $\beta$ is a good matching and $k \geq 0$ (to see the existence of such an expression, simply consider the arc $(i,j) \in C(\alpha)$ with maximal distance between its endpoints; then $k$ corresponds to one of the endpoints). Using Lemma \ref{class_easy}, $[\Psi_{\beta}] = \sum_{I \in T(\beta)} \text{sgn}(I) [\Lambda_{I(1), \cdots, I(m+2n)}]$, and we may compute as follows. Let $\Lambda_{-\lambda} = \bigotimes_{1 \leq s \leq m+2n} \Lambda_{s}^{-\lambda_s}$ and let $T_{-\lambda_1, \cdots, -\lambda_{m+2n}}$ be the operator on $K^0({\rm Coh}(Y_{m+2n}))$ induced by tensoring with the line bundle $\Lambda_{-\lambda}$. For the second-last equality, see the proof of Lemma \ref{computation}; for the last equality, we use Lemma \ref{Eulerchar}.

\begin{align*} d_{\alpha}(\tilde{\lambda}) &= \chi(\Psi(r_{m+2n}^k) \Psi_{\beta} \otimes \Lambda_{-\lambda}) = \sum_{I \in T(\beta)} \text{sgn}(I) \chi(\Psi(r_{m+2n}^k) \Lambda_{I(1), \cdots, I(m+2n)} \otimes \Lambda_{-\lambda}) \\ &= \sum_{I \in T(\beta)} \text{sgn}(I) \chi(T_{- \lambda_1, \cdots, - \lambda_{m+2n}} r_{m+2n}^k v_I) \\ &= \sum_{I \in T(\beta)} \text{sgn}(I) \chi(T_{-\lambda_1, \cdots, - \lambda_{m+2n}} \{ v_{I(k+1)} \otimes \cdots \otimes v_{I(m+2n)} \otimes \theta^{-1}(v_{I(1)}) \otimes \cdots \otimes \theta^{-1}(v_{I(k)})\}) \\ &= \sum_{I \in T(\beta)} \text{sgn}(I) \chi(\theta^{-\lambda_1} v_{I(k+1)} \otimes \cdots \otimes \theta^{-\lambda_{m+2n-k}} v_{I(m+2n)} \otimes \theta^{-\lambda_{m+2n-k+1}-1} v_{I(1)} \otimes \cdots \otimes \theta^{-\lambda_{m+2n}-1} v_{I(k)}) \\ &= \sum_{I \in T(\beta)} \text{sgn}(I) \prod_{i \in I, i \leq k} (\tilde{\lambda}_{m+2n-k+i} + 1) \prod_{i \in I, i > k} \tilde{\lambda}_{i-k} \end{align*} 

We can simplify this further via the following combinatorial observation: 
$$ \sum_{I \in T(\beta)} \text{sgn}(I) \prod_{i \in I} \mu_i = \prod_{(i,j) \in C(\beta)} (\mu_i - \mu_j) $$ 
We finish the proof using this, and the following observation: if $(i,j) \in C(\beta), k \leq i < j$, then $(i-k, j-k) \in C_0(\alpha)$; if $(i,j) \in C(\beta), i < j \leq k$, then $(m+2n-k+i, m+2n-k+j) \in C_0(\alpha)$; if $(i,j) \in C(\beta), i \leq k < j$, then $(j-k, m+2n-k+i) \in C_1(\alpha)$, we get the desired formula: 
\begin{align*} d_{\alpha}(\tilde{\lambda}) &= \prod_{(i,j) \in C(\beta), i,j \leq k} \{ (\tilde{\lambda}_{m+2n-k+i}+1) - (\tilde{\lambda}_{m+2n-k+j}+1)  \} \\  &\prod_{(i,j) \in C(\beta), i,j>k} (\tilde{\lambda}_{i-k} - \tilde{\lambda}_{j-k}) \prod_{(i,j) \in C(\beta), i<k \leq j} (\tilde{\lambda}_{m+2n-k+i}+1 - \tilde{\lambda}_{j-k}) \\ &= \prod_{(i,j) \in C_0(\alpha)} (\tilde{\lambda}_i - \tilde{\lambda}_j) \prod_{(i,j) \in C_1(\alpha)} (1 - \tilde{\lambda}_i + \tilde{\lambda}_j) \end{align*} \end{proof}


\subsection{Dimension formulae}
Recall the definition of $M_{\alpha}$ from Section \ref{[AN]}; and the Main Theorem:

\begin{theorem}[Main Theorem] \label{main} Suppose $\tilde{\mu} \in \mathfrak{h}^*_{\mathbb{Z}}$ satisfies $\langle \mu+\rho, \check{\alpha}_0 \rangle < p$.
\begin{align*} \text{dim}(T_{0 \rightarrow \mu} M_{\alpha}) &= p^{\frac{1}{2}(m+2n)(m+2n-1) - n} \prod_{(i, j) \in C_0(\alpha)} (\tilde{\mu}_i - \tilde{\mu}_j - i +j) \prod_{(i, j) \in C_1(\alpha)} (p-\tilde{\mu}_i + \tilde{\mu}_j + i - j) \end{align*} \end{theorem} \begin{remark} Recall that: $$ \rho = \frac{n-1}{2} e_1 + \frac{n-3}{2} e_2 + \cdots + \frac{1-n}{2} e_n $$ Since $\langle \mu+\rho, \check{\alpha}_0 \rangle \leq p$, note that: \begin{align*} p-\tilde{\mu}_i + \tilde{\mu}_j + i-j \geq p-\tilde{\mu}_1+\tilde{\mu}_n+1-n = p - \langle \mu+\rho, \check{\alpha}_0 \rangle > 0 \end{align*} Hence all terms in the above product are positive. \end{remark} \begin{proof} 
In the following chain of equalities, we use Proposition \ref{dimension} and Lemma \ref{6.2.6} in the first line, and Lemma \ref{Eulertensoring} in the second line.
\begin{align*} \text{dim}(T_{0 \rightarrow \mu} M_{\alpha}) &= \chi(\text{Fr}^* \Psi_{\alpha} \otimes \mathcal{O}(\mu + p \rho)) = p^{\frac{1}{2}(m+2n)(m+2n-1)} d_{\alpha}(\frac{\tilde{\mu} + \rho}{p}) \\ &= p^{\frac{1}{2}(m+2n)(m+2n+1)} \prod_{(i,j) \in C_0(\alpha)} (\frac{\tilde{\mu}_i + \rho_i}{p} - \frac{\tilde{\mu}_j + \rho_j}{p}) \prod_{(i,j) \in C_1(\alpha)} (1 - \frac{\tilde{\mu}_i + \rho_i}{p} + \frac{\tilde{\mu}_j + \rho_j}{p}) \\ &= p^{\frac{1}{2}(m+2n)(m+2n-1) -n} \prod_{(i, j) \in C_0(\alpha)} (\tilde{\mu}_i - \tilde{\mu}_j - i +j) \prod_{(i, j) \in C_1(\alpha)} (p-\tilde{\mu}_i + \tilde{\mu}_j + i - j) \end{align*} \end{proof} 

\subsection{The sub-regular case} \label{subregular} As an application, in the $n=1$ case of Theorem \ref{main} we recover dimension formulae due to Jantzen in the case where the nilpotent is a sub-regular nilpotent. Here $\mathfrak{g}=\mathfrak{sl}_{m+2}$ and the nilpotent has Jordan type $(m+1, 1)$. 

Jantzen shows that there is one irreducible module $L_{\alpha} \in \text{Mod}_{e, 0}(\textbf{U}\mathfrak{sl}_{m+2})$ for each $\alpha \in \Delta^+ \cup \{ \theta \}$; here $\theta$ is the long root, and $\Delta^+ = \{ \epsilon_1 - \epsilon_2, \epsilon_2 - \epsilon_3, \cdots, \epsilon_{m+1} - \epsilon_{m+2} \}$ is the set of simple roots. Jantzen's dimensions formulae (see page $30$ of \cite{jantzen2}) are as follows. \begin{align*} \text{dim}(T_{0 \rightarrow \mu} L_{\alpha}) &=  p^{\frac{(m+1)(m+2)}{2} - 1} \langle \mu + \rho, \alpha \rangle \text{ if } \alpha \in \Delta^+ \\ \text{dim}(T_{0 \rightarrow \mu} L_{\theta}) &= p^{\frac{(m+1)(m+2)}{2} - 1}  (p - \langle \mu + \rho,\check{\theta} \rangle) \end{align*} 

Using our parametrization, the irreducibles are as follows. For $1 \leq i \leq m+1$, let $M_i$ be the simple object indexed by the unique crossingless matching in $\text{Cross}(m,1)$ with a cup joining $i$ and $i+1$. Let $M_{m+2}$ be the simple object indexed by the unique crossingless matching in $\text{Cross}(m,1)$ with a cup joining $m+2$ and $1$. Using Theorem \ref{main}: 
\begin{align*} \text{If } 1 \leq i \leq m+1: \text{dim}(T_{0 \rightarrow \mu} M_i) &= p^{\frac{(m+1)(m+2)}{2} - 1} (\tilde{\mu}_i - i - \tilde{\mu}_{i+1} + i +1)  = p^{\frac{(m+1)(m+2)}{2} - 1} (\tilde{\mu}_i - \tilde{\mu}_{i+1} + 1) \\  \text{dim}(T_{0 \rightarrow \mu} M_{m+2}) &= p^{\frac{(m+1)(m+2)}{2} - 1} (p - \tilde{\mu}_1 + \tilde{\mu}_{m+2} - m - 1) \end{align*}
It follows that $L_{\epsilon_i - \epsilon_{i+1}} \simeq M_i$ and $L_{\theta} \simeq M_{m+2}$ (since $\langle \rho, \check{\theta} \rangle = m+1$, $\langle \mu, \check{\theta} \rangle = \tilde{\mu}_1 - \tilde{\mu}_{m+2}$). This completes the verification of Jantzen's formulae in type A. 

\subsection{The case when $n=2$}

In this subsection we let $n=2$ and let $m$ be an arbitrary positive integer, and we look at the dimension formulas which the Main Theorem gives us in this case. We have that $\alpha$ is determined up to isotopy by the set $C(\alpha)$. So we will compute the dimensions of the modules in the Main Theorem for all possible $C(\alpha)$.

Case 1. $C(\alpha)=\{(i,i+1),(j,j+1): 1 \leq i, i+1<j, j+1 \leq m+4\}$. Then $C_0(\alpha)=C(\alpha)$ and $C_1(\alpha)=\emptyset$. Then we have:
$$\text{dim}(T_{0 \rightarrow \mu} M_{\alpha})   = p^{\frac{(m+3)(m+4)}{2} - 2} (\tilde{\mu}_i - \tilde{\mu}_{i+1} + 1)(\tilde{\mu}_j - \tilde{\mu}_{j+1} + 1).$$

Case 2. $C(\alpha)=\{(i,i+1),(1,m+4): 1 < i < m+3\}$. Then $C_0(\alpha)=\{(i,i+1)\}$ and $C_1(\alpha)=\{(1,m+4)\}$. Then we have:
$$\text{dim}(T_{0 \rightarrow \mu} M_{\alpha})   = p^{\frac{(m+3)(m+4)}{2} - 2} (\tilde{\mu}_i - \tilde{\mu}_{i+1} + 1)(p - \tilde{\mu}_1 + \tilde{\mu}_{m+4} - m - 3).$$

Case 3. $C(\alpha)=\{(i,i+3),(i+1,i+2): 1 \leq i \leq m+1\}$. Then $C_0(\alpha)=C(\alpha)$ and $C_1(\alpha)=\emptyset$. Then we have:
$$\text{dim}(T_{0 \rightarrow \mu} M_{\alpha})   = p^{\frac{(m+3)(m+4)}{2} - 2} (\tilde{\mu}_i - \tilde{\mu}_{i+3} + 3)(\tilde{\mu}_{i+1} - \tilde{\mu}_{i+2} + 1).$$

Case 4. $C(\alpha)=\{(3,m+4),(1,2)\}$. Then $C_0(\alpha)=\{(1,2)\}$ and $C_1(\alpha)=\{(3,m+4)\}$. Then we have:
$$\text{dim}(T_{0 \rightarrow \mu} M_{\alpha})   = p^{\frac{(m+3)(m+4)}{2} - 2} (\tilde{\mu}_1 - \tilde{\mu}_2 + 1)(p - \tilde{\mu}_3 + \tilde{\mu}_{m+4} - m - 1).$$

Case 5. $C(\alpha)=\{(1,m+2),(m+3,m+4)\}$. Then $C_0(\alpha)=\{(m+3,m+4)\}$ and $C_1(\alpha)=\{(1,m+2)\}$. Then we have:
$$\text{dim}(T_{0 \rightarrow \mu} M_{\alpha})   = p^{\frac{(m+3)(m+4)}{2} - 2} (\tilde{\mu}_{m+3} - \tilde{\mu}_{m+4} + 1)(p - \tilde{\mu}_1 + \tilde{\mu}_{m+2} - m - 1).$$

Case 6. $C(\alpha)=\{(1,m+4),(2,m+3)\}$. Then $C_0(\alpha)=\emptyset$ and $C_1(\alpha)=C(\alpha)$. Then we have:
$$\text{dim}(T_{0 \rightarrow \mu} M_{\alpha})   = p^{\frac{(m+3)(m+4)}{2} - 2} (p - \tilde{\mu}_1 + \tilde{\mu}_{m+4} - m - 3)(p - \tilde{\mu}_2 + \tilde{\mu}_{m+3} - m - 1).$$

\subsection{Kac-Weisfeiler theory}

In this subsection, we return to the setting of Section \ref{modularrep}; let $\mathfrak{g}$ be a semisimple Lie algebra and $e \in \mathfrak{g}^*$ an arbitrary nilpotent. Assume in addition that $e$ is in standard Levi form. Recall that the irreducible modules in $\text{Mod}^{\text{fg}}_{e, \lambda}(U \mathfrak{g})$ are of the form $L_e(\mu)$. 

\begin{definition} Given a nilpotent $e \in \mathfrak{g}^*$, let $d(e)$ be the co-dimension of the Springer fiber $\mathcal{B}_e$ in the flag variety $\mathcal{B}$. \end{definition}

The Kac-Weisfeiler conjecture (now a theorem of Premet, \cite{pr}), was originally stated for arbitrary prime $p$; we are more interested in the case where $p \gg 0$. In this case, in \cite{bmr} Bezrukavnikov-Mirkovic-Rumynin gave another proof using localization techniques. 

\textbf{Kac-Weisfeiler conjecture.} The dimension of any irreducible module in $\text{Mod}_{e, \lambda}^{fg}(U \mathfrak{g})$ is divisible by $p^{d(e)}$. 

Below we verify that our results are consistent with the Kac-Weisfeiler conjecture when $e$ is a two-row nilpotent and $p \gg 0$. 

\begin{proof}[Proof (when $e$ is a two-row nilpotent in type A)] It suffices to show that the dimension of any irreducible $L_e(\mu)$ is divisible by $p^{d(e)}$. Using the formula for dimensions of Springer fibers in type A (see Spaltenstein, \cite{spa}), it follows that $\text{dim}(\mathcal{B}_{e}) = n$.
 
$$\text{codim}_{\mathcal{B}}(\mathcal{B}_{e}) = \text{dim}(\mathcal{B}) - \text{dim}(\mathcal{B}_{e}) =  \frac{1}{2}(m+2n)(m+2n-1) -n$$
 
Hence Theorem \ref{main} shows that the Kac-Weisfeiler conjecture is true in this case. \end{proof} 

\textbf{Humphreys' conjecture (see \cite{hum}).} For a fixed nilpotent $e \in \mathfrak{g}^*$, there is at least one module $L_e(\mu)$ with dimension equal to $p^{d(e)}$. 

This conjecture is now known - see Theorem 2.2 in \cite{pr2}; see also recent work of Goodwin and Topley, \cite{gt}, for an extension of this result counting the number of such modules. Below we verify that this is consistent with our Theorem \ref{main}. 

\begin{proof}[Proof (when $e$ is a two-row nilpotent in type A)] Choose any $\alpha \in \text{Cross}(m,n)$ with the following property: if $1 \leq i,j \leq m+2n$, and $i$ and $j$ are connected in $\alpha$, then $|i-j|=1$. The existence of such an $\alpha$ is clear (when $m=0$, there is exactly one such $\alpha$). Using Theorem \ref{main}, any such $\alpha$ satisfies $\text{dim}(M_{\alpha}) = p^{d(e)}$. \end{proof}

\section{Multiplicities of simples in baby Vermas} \label{multip}

In the notation of Section \ref{[AN]}, let $\Delta$ be a baby Verma module in $\text{Mod}^{\text{fg}}_{e, 0}(U \mathfrak{g})$. Recall from Section \ref{modularrep} that given a simple module $M_{\alpha}$, the composition multiplicity $[\Delta: M_{\alpha}]$ does not depend on the choice of baby Verma, i.e. it only depends on $\alpha$. Denote this quantity by $m(\alpha)$. 

\begin{remark} These multiplicities are preserved under the action of translation functors, so the theorem below would remain true if $\Delta$ (resp. $M_{\alpha}$) was replaced by a baby Verma (resp. corresponding simple) in $\text{Mod}^{\text{fg}}_{e, \lambda}(U \mathfrak{g})$, where $\lambda$ is a regular Harish-Chandra character. \end{remark} 

\begin{definition} Given $\alpha \in \text{Cross}(m,n)$, given two distinct arcs $(i,j)$ and $(k,l)$, say that $(k,l) < (i,j)$ if the arc $(i,j)$ separates the arc $(k,l)$ from the inner circle. Denote by $c(i,j)$ the number of arcs $(k,l) \in C(\alpha)$, such that either $(k,l) = (i,j)$, or $(k,l) < (i,j)$. \end{definition}

\begin{theorem} \label{mult} $$m(\alpha) = \frac{n!}{\displaystyle \prod_{(i,j) \in C(\alpha)} c(i,j)}$$ \end{theorem} 

\begin{example} In Example \ref{examp}, $c(1,2) = 1, c(3,6)=2$ and $m(\alpha)=1$. \end{example}

\begin{definition} \label{label} A labelling of the $n$ arcs in $C(\alpha)$ with the numbers $\{ 1, \cdots, n \}$ is called ``good'' if for every pair of arcs $(i,j), (k,l)$ with $(k,l) < (i,j)$, then the label attached to $(k,l)$ is less than the label attached to $(i,j)$. Let $\underline{\text{Cross}}(m,n)$ denote the set of all crossingless matchings $\alpha \in \text{Cross}(m,n)$ equipped with a good labelling $L$, and let $n(\alpha)$ be the number of good labellings of $\alpha$. \end{definition}

Recall from Section \ref{embedding} that we identified $K^0(\text{Mod}^{\text{fg}}_{e, 0}(U \mathfrak{g}))$ with a subspace of $(\mathbb{C}^2)^{\otimes m+2n}$. Under this identification, $[M_{\alpha}] = \underline{\alpha}$. 

\begin{lemma} \label{sky} $$ [\Delta] = v_0 \otimes \cdots \otimes v_0 $$ \end{lemma} \begin{proof} In Section $5.3.3$ of \cite{bmr}, Bezrukavnikov, Mirkovic and Rumynin show that the baby Verma modules correspond to skyscraper sheaves under their localization equivalence. Using Lemma \ref{base}, the class of a skyscraper sheaf is $v_0 \otimes \cdots \otimes v_0$, as required. \end{proof} 

\begin{proposition} \label{classes} $$\sum_{\alpha \in \text{Cross}(m,n)} n(\alpha) \underline{\alpha} = v_0 \otimes \cdots \otimes v_0 $$ \end{proposition} \begin{proof} Proof by induction on $n$. Consider the operator $$F_{m,n}: (\mathbb{C}^2)^{\otimes m+2n} \rightarrow (\mathbb{C}^2)^{\otimes m+2n+2}, \qquad F_{m,n}(v) = (1+\overline{r_{m+2n}}+\cdots+\overline{r_{m+2n}}^{m+2n-1}) g_{m+2n}^1  (v)$$ For the induction step, it suffices to prove the following two equations: \begin{align*} F_{m,n}(\underbrace{v_0 \otimes \cdots \otimes v_0}_{m+2n}) &= \underbrace{v_0 \otimes \cdots \otimes v_0}_{m+2n+2} \\ F_{m,n}(\sum_{\alpha \in \text{Cross}(m,n)} n(\alpha) \underline{\alpha}) &= \sum_{\beta \in \text{Cross}(m,n+1)} n(\beta) \underline{\beta} \end{align*}
The first equation is a computation using Proposition \ref{tanglecat} and Lemma \ref{computation}:  
\begin{align*} F_{m,n}(\underbrace{v_0 \otimes \cdots \otimes v_0}_{m+2n}) &= (1+\overline{r_{m+2n}}+\cdots+\overline{r_{m+2n}}^{m+2n-1})  (v_1 \otimes v_0 \otimes v_0 \otimes \cdots \otimes v_0 - v_0 \otimes v_1 \otimes v_0 \otimes \cdots \otimes v_0) \\ &= v_0 \otimes \theta^{-1}(v_1) \otimes v_0 \otimes \cdots \otimes v_0 - v_0 \otimes v_1 \otimes v_0 \otimes \cdots \otimes v_0 \\ &= \underbrace{v_0 \otimes \cdots \otimes v_0}_{m+2n+2}  \end{align*} 
In the second equality above we use the following explicit formula and cancelled like terms. $$\overline{r_{m+2n}}^k v_I = v_{I(k+1)} \otimes \cdots \otimes v_{I(m+2n)} \otimes \theta^{-1}(v_{I(1)}) \otimes \cdots \otimes \theta^{-1}(v_{I(k)})$$
   The second equation, re-expressed below, is true since every labelled matching $(\beta, L') \in \underline{\text{Cross}}(m,n+1)$ can be uniquely expressed as $(\beta, L') = r_{m+2n+2}^k g_{m+2n}^1 (\alpha, L)$ for some $(\alpha, L) \in \underline{\text{Cross}}(m,n)$ and an integer $0 \leq k \leq m+2n+1$. Here $r_{m+2n+2}^k$ acts naturally on the set of labelled diagrams, and by $g_{m+2n}^1 (\alpha, L)$ we mean the labelled diagram obtained by labelling the arc $(1,2)$ with a $1$, and increasing all the other labels by $1$.
$$ F_{m,n}(\sum_{(\alpha, L) \in \underline{\text{Cross}}(m,n)} \alpha) = \sum_{(\beta, L') \in \underline{\text{Cross}}(m,n+1)} \beta $$
\end{proof}

\begin{lemma} \label{count} $$ n(\alpha) = \frac{n!}{\displaystyle \prod_{(i,j) \in C(\alpha)} c(i,j)} $$ \end{lemma} This Lemma is a special case of the following, applied to the the poset whose vertices are elements of $C(\alpha)$, with the order relation specified in Definition \ref{label}. The below Lemma is well-known, but we include a proof for completeness. \begin{lemma} Let $P$ be a poset whose Hasse diagram is a tree such that every connected component has a maximal element. Let $V$ be the set of vertices of $P$. If $P$ has $n$ vertices, a labelling of its vertices is a bijection from $V$ to $\{ 1, \cdots, n \}$. A labelling is called ``good'' if for every two vertices $a,b$ with $a \leq b$, then $l(a) \leq l(b)$. Given a vertex $a$, let $$c(a) = \# \{ b \in V, b \leq a \} $$ Then the number of good labellings is $$ f(P) = \frac{n!}{\displaystyle \prod_{1 \leq i \leq n} c(i)} $$ \end{lemma} \begin{proof} Proof by induction on the number of vertices. Suppose that the Hasse diagram has $k$ connected components $R_1, \cdots R_k$ of sizes $r_1, \cdots, r_k$. If $k>1$, there are $\frac{n!}{r_1! r_2! \cdots r_k!}$ ways of choosing the sets $l(R_1), \cdots, l(R_k)$, and once this is done there are $f(R_1) \cdots f(R_k)$ ways of choosing the label function (using the induction hypothesis applied to those components). The conclusion follows since: $$\frac{n!}{r_1! r_2! \cdots r_k!} f(R_1) \cdots f(R_k) = f(P)$$ If $k=1$, then the unique maximal element has label $n$. Let $P'$ be the poset obtained by removing this element; then $f(P) = f(P')$, and the number of good labellings of $P$ is equal to the number of good labellings of $P'$. The conclusion follows using the inductive hypothesis, completing the proof. \end{proof} 

\begin{proof}[Proof of Theorem \ref{mult}] It suffices to prove the equality on the level of Grothendieck groups, i.e. that $$ [\Delta] = \sum_{\alpha \in \text{Cross}(m,n)} \frac{n!}{\displaystyle \prod_{(i,j) \in C(\alpha)} c(i,j)} [M_{\alpha}] $$ This follows from Lemma \ref{sky}, Proposition \ref{classes}, and Lemma \ref{count}. \end{proof}

\bibliographystyle{plain}

\end{document}